\documentclass{article}
\usepackage[left=2.5cm,right=2.5cm,top=3cm,bottom=2.5cm]{geometry}
\usepackage{amsmath,amsthm,amssymb}
\usepackage{bm}

\newtheorem{theorem}{Theorem}[section]
\newtheorem{corollary}[theorem]{Corollary}
\newtheorem{lemma}[theorem]{Lemma}

\newcommand{\pt}{\partial}
\DeclareMathOperator{\dive}{div}

\newcounter{mystep}
\newcommand{\myemph}{\noindent\it \stepcounter{mystep} {\rm (\themystep)} }

\date{}

\begin{document}


\title{The well-posedness of three-dimensional Navier-Stokes and magnetohydrodynamic equations with partial fractional diffusion} 

\author{Qibo Ma$^*$ and Li Li\footnote{School of Mathematics and Statistics, Ningbo University, 818 Fenghua Road, Ningbo 315211, China. Email: lili2@nbu.edu.cn}}




\maketitle

\begin{abstract}
	It is well-known that if one replaces standard velocity and magnetic diffusion by $(-\Delta)^\alpha u$ and $(-\Delta)^\beta b$ respectively, the magnetohydrodynamic equations are well-posed for $\alpha\ge\frac{5}{4}$ and $\alpha + \beta \ge \frac{5}{2}$. This paper considers the 3D Navier-Stokes and magnetohydrodynamic equations with partial fractional hyper-diffusion. It is proved that when each component of the velocity and magnetic field lacks diffusion along some direction, the existence and conditional uniqueness of the solution still hold. This paper extends the previous results in (Yang, Jiu and Wu J. Differential Equations 266(1): 630-652, 2019) to a more general case. 

	\medskip
	\noindent
	{\it Keywords:}
		Navier-Stokes equations, magnetohydrodynamic equations, partial diffusion, well-posedness, fractional diffusion

\end{abstract}

\section{Introduction}

In this paper, we first consider the well-posedness of the following incompressible Navier-Stokes equation with partial fractional hyper-diffusion, 
\begin{equation}\label{eq:general:NS}
	\left\{ 
	\begin{array}{l}
		\partial _t u_1 + \boldsymbol{u} \cdot \nabla u_1 + \nu (\Lambda_1 ^{5/2} u_1 + \Lambda _2^{5/2} u_1 + \Lambda _3^{5/2} u_1 - \Lambda _{i_1}^{5/2} u_1) + \partial _1 p = 0 ,\\
		\partial _t u_2 + \boldsymbol{u} \cdot \nabla u_2 + \nu (\Lambda_1 ^{5/2} u_2 + \Lambda _2^{5/2} u_2 + \Lambda _3^{5/2} u_2 - \Lambda _{i_2}^{5/2} u_2) + \partial _2 p = 0,\\
		\partial _t u_3 + \boldsymbol{u} \cdot \nabla u_3 + \nu (\Lambda_1 ^{5/2} u_3 + \Lambda _2^{5/2} u_3 + \Lambda _3^{5/2} u_3 - \Lambda _{i_3}^{5/2} u_3) + \partial _3 p = 0,\\
		\dive \boldsymbol{u} = 0,\\
		\boldsymbol{u}( x,0 ) = \boldsymbol{u}_0( x ), 
	\end{array} 
	\right.
\end{equation}
where $ {\bm u} = (u_1,u_2,u_3) $ represents the velocity field, $p$ is the pressure and $\nu$ denotes the viscosity, $i_1,i_2,i_3\in\{1,2,3\}$, the directional fractional operator $ \Lambda _k^\gamma $ is defined via Fourier transform 
\begin{equation}\label{eq:Lambda:k}
	\widehat {\Lambda _k^\gamma f( \xi )} := | \xi _k |^\gamma \widehat {f( \xi )},\quad \gamma>0, k=1, 2, 3. 
\end{equation}
The above equation implies that the velocity component $u_k$ lacks of diffusion along the $x_{i_k}$ direction for $k=1,2,3$. 

We further study the well-posedness of the three dimensional incompressible magnetohydrodynamic (MHD) equations with partial hyper-diffusion,
\begin{equation}\label{eq:general:MHD}
	\left\{ 
	\begin{array}{l}
		\partial _t u_1 + \boldsymbol{u} \cdot \nabla u_1 + \nu (\Lambda_1 ^{5/2} u_1 + \Lambda _2^{5/2} u_1 + \Lambda _3^{5/2} u_1 - \Lambda _{i_1}^{5/2} u_1) + \partial _1 p = {\bm b}\cdot\nabla b_1, \\
		\partial _t u_2 + \boldsymbol{u} \cdot \nabla u_2 + \nu (\Lambda_1 ^{5/2} u_2 + \Lambda _2^{5/2} u_2 + \Lambda _3^{5/2} u_2 - \Lambda _{i_2}^{5/2} u_2) + \partial _2 p = {\bm b}\cdot\nabla b_2, \\
		\partial _t u_3 + \boldsymbol{u} \cdot \nabla u_3 + \nu (\Lambda_1 ^{5/2} u_3 + \Lambda _2^{5/2} u_3 + \Lambda _3^{5/2} u_3 - \Lambda _{i_3}^{5/2} u_3) + \partial _3 p = {\bm b}\cdot\nabla b_3, \\
		\pt_t {\bm b} + {\bm u} \cdot \nabla {\bm b} + \eta ( \Lambda_1 ^{5/2} {\bm b} + \Lambda_2^{5/2} {\bm b} + \Lambda_3 ^{5/2} {\bm b} - \Lambda_{j_0} ^{5/2} {\bm b} ) = {\bm b} \cdot \nabla {\bm u}, \\
		\dive \boldsymbol{u} = 0, \qquad \dive {\bm b} = 0,\\
		\boldsymbol{u}( x,0 ) = \boldsymbol{u}_0( x ), \qquad {\bm b}( x,0 ) = {\bm b}_0( x ), 
	\end{array} 
	\right.
\end{equation}
where $ {\bm u} $ and $ {\bm b} $ represent the velocity and magnetic field, respectively, $\nu$ and $\eta$ denote the viscosity and magnetic diffusion coefficients, respectively, $i_k,j_0 \in\{1,2,3\}$. The above equation implies that $u_k$ lacks of diffusion along the $x_{i_k}$ direction, while the magnetic component $b_k$, $k=1,2,3$, lacks of diffusion along the $x_{j_0}$ direction. It should be noted that, in this generalized MHD system, the divergence-free condition for the magnetic field are kept only if each component of ${\bm b}$ lacks of the diffusion in the same direction. 

The Navier-Stokes equations describe the motion of viscous fluids, and the well-posedness problem of the three-dimensional incompressible equations remains an open challenge. In 1934, Leray \cite{Leray1934} proved the existence of global weak solutions for the incompressible Navier-Stokes equations. Thus far, the best regularity result is attributed to Caffarelli, Kohn, and Nirenberg's work \cite{CKN1982}, showing that singular sets of suitable weak solutions have one-dimensional Hausdorff measure zero. Later, Lin \cite{Lin1998} proved the same result using compactness methods. Vasseur \cite{Vasseur2007} also demonstrated the same result by utilizing the techniques introduced by De Giorgi for elliptic equations. In 2019, Paicu and Zhang \cite{PZ2019} proved the existence of global strong solutions for the three-dimensional incompressible Navier-Stokes equations when the initial data satisfies the anisotropic small condition. This condition takes into account the different effects of horizontal and vertical viscosities, see also \cite{Zhang2009,LZ2020}. In 2021, Liu and Zhang \cite{LZ2021} proved that there exists a unique global strong solution for the three-dimensional anisotropic incompressible Navier-Stokes equations with strong dissipation in the vertical direction. Furthermore, the norm of the vertical component of the velocity field can be controlled by the norm of the corresponding component of the initial data. 

Due to insufficient dissipation to control the nonlinearity of the equations, the standard Navier-Stokes equations are considered to be supercritical. 
Katz and Pavlovi\'{c} \cite{KP2002} considered the following generalized Navier-Stokes equations with fractional hyper-diffusion 
\[
\left\{
\begin{array}{l}
	\partial _t \bm {u} + {\bm u} \cdot \nabla {\bm u} = - \nabla p - \nu (-\Delta)^\alpha {\bm u}, \qquad x \in \mathbb{R}^3, t>0,\\
	\dive {\bm u} = 0,\\
	{\bm u}( x,0 ) = {\bm u}_0(x), 
\end{array} 
\right.
\]
where $\alpha \ge \frac{5}{4}$ and $(-\Delta)^\alpha$ is the fractional Laplacian operator 
\[
	\widehat {(-\Delta)^\alpha f(\xi)} := |\xi|^{2\alpha }\widehat {f( \xi )}.
\]
It was proved that any smooth initial value ${\bm u}_0$ with finite energy will have a unique global-in-time solution. See also \cite{Wu2003,Lions1969}. 
Later, Tao improved the above results in \cite{Tao2009}. He studied the Navier-Stokes equation with the following logarithmic form of supercritical dissipative term, 
\[
\left\{ 
\begin{array}{l}
	\partial _t {\bm u} + {\bm u} \cdot \nabla {\bm u} = - \nabla p - \nu \frac{(-\Delta)^\frac{5}{4}}{\log_{}^{\frac{1}{2}}(2-\Delta)} {\bm u}, \qquad x \in \mathbb{R}^3, t>0,\\
	\dive {\bm u} = 0,\\
	{\bm u}(x,0) = {\bm u}_0( x ),
\end{array} 
\right.
\]
and found that the above problem has a global smooth solution for any compactly supported smooth initial data. Barbato, Morandin and Romito \cite{BMR2014} improved the research results of Tao by introducing weaker conditions on the logarithmic weak dissipative term. The result was further improved by Yamazaki to MHD equaion in \cite{Yamazaki2018}. 
Yang, Jiu and Wu in \cite{YJW2019} reduced some excessive dissipation. They proved the global regularity of the following 3D Navier-Stokes equations with partial dissipation,
\[
\left\{ 
\begin{array}{l}
	\partial _t u_1 + {\bm u} \cdot \nabla u_1 = - \partial _1 p - \nu \Lambda _1^{5/2} u_1 - \nu \Lambda _2^{5/2} u_1,\\
	\partial _t u_2 + {\bm u} \cdot \nabla u_2 = - \partial _2 p - \nu \Lambda _2^{5/2} u_2 - \nu \Lambda _3^{5/2} u_2,\\
	\partial _t u_3 + {\bm u} \cdot \nabla u_3 = - \partial _3 p - \nu \Lambda _1^{5/2} u_3 - \nu \Lambda _3^{5/2} u_3,\\
	\dive {\bm u} = 0,\\
	{\bm u}( x,0 ) = {\bm u}_0( x ),
\end{array} 
\right.
\]
where $ \Lambda _k^\gamma $ denotes the directional fractional operators defined in \eqref{eq:Lambda:k}. 
However, when all equations loss dissipation in the same direction, for example, the well-posedness of 
\[
\left\{ 
\begin{array}{l}
	\partial _t u_1 + {\bm u} \cdot \nabla u_1 = - \partial _1 p - \nu \Lambda _1^{5/2} u_1 - \nu \Lambda _2^{5/2} u_1,\\
	\partial _t u_2 + {\bm u} \cdot \nabla u_2 = - \partial _2 p - \nu \Lambda _1^{5/2} u_2 - \nu \Lambda _2^{5/2} u_2,\\
	\partial _t u_3 + {\bm u} \cdot \nabla u_3 = - \partial _3 p - \nu \Lambda _1^{5/2} u_3 - \nu \Lambda _2^{5/2} u_3,\\
	\dive {\bm u} = 0,\\
	{\bm u}( x,0 ) = {\bm u}_0( x ),
\end{array} 
\right.
\]
can not be obtained through the method in \cite{YJW2019}. 

In this paper, we first study the existence and uniqueness of equation \eqref{eq:general:NS}. 
Let 
\begin{equation}\label{eq:I}
	\mathcal{I} := \big\{(i_1, i_2, i_3) \mid i_k \in\{1, 2, 3\}, i_1=i_2=3 \text{ or } i_2=i_3=1 \text{ or } i_1=i_3=2 \big\}
\end{equation}
be the set of all "bad" indices for ${\bm u}$. 
Our main results on Naiver-Stokes equations are the following theorems. 
\begin{theorem}[Existence]\label{thm:existence:NS}
	If the indices $(i_1, i_2, i_3) \notin \mathcal{I}$ and $ \boldsymbol{u}_0 \in H^1( \mathbb{R}^3 ) $, then there exists a global solution ${\bm u}$ of \eqref{eq:general:NS} satisfying
	\begin{equation}\label{eq:thm:exis}
		\begin{array}{l}
			{\bm u} \in L^\infty \big( (0,\infty );H^1\big),\\
			\Lambda_j^{\frac{5}{4}} u_1, \nabla \Lambda_j ^\frac{5}{4} u_1 \in L^2 \left( (0,\infty) ;L^2\right), \quad j=1,2,3, \quad j\not= i_1,\\
			\Lambda_j^\frac{5}{4} u_2, \nabla \Lambda_j^\frac{5}{4} u_2 \in L^2\left( (0,\infty) ;L^2 \right), \quad j=1,2,3, \quad j\not= i_2,\\
			\Lambda_j^\frac{5}{4} u_3, \nabla \Lambda_j^\frac{5}{4} u_3 \in L^2\left( (0,\infty) ;L^2 \right), \quad j=1,2,3,\quad j\not= i_3.
		\end{array}			
	\end{equation}
	In addition, the $H^1$-norm of ${\bm u}$ is bounded uniformly in time.
\end{theorem}
To ensure uniqueness, one needs to make additional assumptions for certain indices that one of the solutions should has better regularity in $H^1$ space. 
\begin{theorem}[Uniqueness]\label{thm:uniqueness:NS}
	Suppose that $\bm{u}^{(1)}$ and $\bm{u}^{(2)}$ are two solutions of \eqref{eq:general:NS} on $(0,T)$ for some $(i_1, i_2, i_3)$ $\notin \mathcal{I}$, and satisfy \eqref{eq:thm:exis} on $(0,T)$. 
	We further assume that 
	\begin{equation}\label{eq:thm:uniq}
		\Lambda_k^\frac{9}{4}u_l^{(2)}\in L^2(0,T;L^2), 
	\end{equation}
	if there exist $k,l\in\{1,2,3\}$, $k\not=l$, such that $i_k=i_l=k$. 
	Then $\bm{u}^{(1)} \equiv \bm{u}^{(2)}$ in $\mathbb{R}^3\times(0,T)$.
\end{theorem}

\vspace{10pt}

The second part of this article focuses on the well-posedness of the MHD equations with partial hyper-diffusion. 
The MHD equations combine the magnetic field with the Navier-Stokes equations in hydrodynamics to describe the behavior of charged fluids, such as plasma. These equations are essential for studying various phenomena in geophysics, astrophysics, cosmology, and engineering (see \cite{Biskamp1993,PF2000,Davidson2001}). In 1972, Duvaut and Lions \cite{DL1972} proved that there is at least one global weak solution for the MHD equations and it is locally well-posed. In 2003, Wu \cite{Wu2003} proved the global regularity in the critical case of weak solutions to the generalized MHD equations. Cao and Wu \cite{CW2011} proved the global regularity results of two-dimensional incompressible MHD equations with partial dissipation and magnetic diffusion in 2011. Later, Wang and Wang in \cite{WW2013} gave the existence of global smooth solutions for the initial value problem of three-dimensional incompressible MHD equations with a mixture of partial dissipation and magnetic diffusion. In 2015, Jiu and Liu \cite{JL2015} studied the Cauchy problem of three-dimensional axisymmetric MHD equations with horizontal dissipation and vertical magnetic diffusion. Under the assumption that $u_\theta$ and $b_r$ are trivial, they obtained a unique globally smooth solution. In 2022, Lin, Wu and Zhu \cite{LWZ2022} dealt with Navier-Stokes nonlinearity when the dissipation of the velocity is only in a single direction and proved that any disturbance near the background magnetic field $(0, 1, 0)$ is globally stable in the Sobolev space $H^3(\mathbb{R}^3)$. The global existence and regularity results for $2\frac{1}{2}$-dimensional Hall-MHD equatios were studied in \cite{RY,JJ2024}. More related results can be seen in \cite{HLL2008,Lei2015,WZ2021} and the references provided therein. 

The three-dimensional incompressible MHD equations with standard Laplace dissipation is considered to be supercritical, and has the same difficulties as the standard Navier-Stokes equations. As is well-known (c.f. \cite{Wu2003,Wu2011,Yamazaki2014}), the generalized MHD equations
\[
\left\{ 
\begin{array}{l}
	\partial _t {\bm u} + {\bm u} \cdot \nabla {\bm u} = - \nabla p - \nu (-\Delta)^\alpha {\bm u}+{\bm b}\cdot \nabla {\bm b}, \\
	\partial _t {\bm b} + {\bm u} \cdot \nabla {\bm b} = - \eta ( - \Delta )^\beta {\bm b}+{\bm b}\cdot \nabla {\bm u}, \\
	\dive {\bm u} = 0,\qquad \dive {\bm b} = 0,\\
\end{array} 
\right.
\]
with $\alpha\ge\frac{5}{4},\alpha+\beta\ge\frac{5}{2}$ always possess global classical solutions corresponding to any sufficiently smooth initial data. 
After removing some excessive diffusion along some directions, the existence and uniqueness of solutions to 
\[
	\left\{ 
	\begin{array}{l}
		\pt_t u_1 + {\bm u} \cdot \nabla u_1 +\nu \Lambda^\frac{5}{2} u_1 + \pt_1 p = {\bm b}\cdot\nabla b_1, \\
		\pt_t u_2 + {\bm u} \cdot \nabla u_2 +\nu ( \Lambda_2^\frac{5}{2} u_2+\Lambda _3^\frac{5}{2} u_2 ) + \pt_2 p = {\bm b}\cdot \nabla b_2, \\
		\pt_t u_3 + {\bm u} \cdot \nabla u_3+\nu ( \Lambda_1 ^\frac{5}{2} u_3 + \Lambda _3^\frac{5}{2} u_3 ) + \pt_3 p = {\bm b}\cdot \nabla b_3, \\
		\pt_t {\bm b} + {\bm u} \cdot \nabla {\bm b} + \eta ( \Lambda_1 ^\frac{5}{2} {\bm b} + \Lambda_2^\frac{5}{2} {\bm b} ) = {\bm b} \cdot \nabla {\bm u}, \\
		\dive {\bm u} = 0, \qquad \dive {\bm b} = 0,\\
		{\bm u}( x,0 ) = {\bm u}_0( x ), \qquad {\bm b}( x,0 ) = {\bm b}_0( x ),
	\end{array} 
	\right.
\]
was proved in \cite{YJW2019M} 
where $ \Lambda := ( - \Delta )^\frac{1}{2} $ denotes the Zygmund operator.
Notice that the above well-posedness result was obtained when the velocity component $u_1$ has dissipation in every direction, $u_2$ lacks dissipation in $x_1$ direction, $u_3$ lacks dissipation in $x_2$ direction, and each component of magnetic field ${\bm b}$ lacks dissipation in $x_3$ direction. 
Inspired by \cite{YJW2019M}, we study the existence and uniqueness of problem \eqref{eq:general:MHD} for indices $i_k, j_0\in\{1,2,3\}$, $k=1,2,3$. 

Denote the set of all "bad" indices as
\begin{equation}\label{eq:J}
	\begin{aligned}
		\mathcal{J} := & \big\{(i_1, i_2, i_3, j_0) \mid i_k, j_0 \in\{1, 2, 3\} \text{ for } k=1,2,3; \\
		& (i_1, i_2, i_3) \in \mathcal{I} \text{ or } j_0 = i_{l} \in\{1,2,3\} \backslash\{l\} \text{ for some } l \in\{1,2,3\}\big\},
	\end{aligned}
\end{equation}
where $\mathcal{I}$ being the "bad" indices for the Navier-Stokes equations defined in \eqref{eq:I}. 
Our main results on MHD equations are the following theorems.
\begin{theorem}[Existence]\label{thm:existence:MHD}
	If $(i_1, i_2, i_3, j_0) \notin \mathcal{J}$ and $(\boldsymbol{u}_0, \boldsymbol{b}_0)\in H^1( \mathbb{R}^3 ) $, then there exists a global solution $({\bm u},{\bm b})$ of \eqref{eq:general:MHD} satisfying
	\begin{equation}\label{eq:reg}
		\begin{array}{l}
			{\bm u} \in L^\infty \big( (0,\infty );H^1\big),\quad {\bm b} \in L^\infty \big( (0,\infty );H^1\big),\\
			\Lambda_k^\frac{5}{4} u_l, \nabla \Lambda_k ^\frac{5}{4} u_l \in L^2 \left( (0,\infty); L^2\right), \quad k,l=1,2,3, \quad k \not= i_l, \\
			\Lambda_k^\frac{5}{4} b_l, \nabla \Lambda_k ^\frac{5}{4} b_l \in L^2 \left( (0,\infty); L^2\right), \quad k,l=1,2,3, \quad k\not= j_0.
		\end{array}
	\end{equation}
	The bound on the norms of ${\bm u }$ and ${\bm b}$ in \eqref{eq:reg} is uniform in time.
\end{theorem}	
Furthermore, we have the following uniqueness result for the MHD equations. Likewise, assumptions must be made for certain indices, requiring one of the solutions to exhibit enhanced regularity in the $H^1$ space. 

\begin{theorem}[Uniqueness]\label{thm:uniqueness:MHD}
	Let $T>0$. Assume that $(i_1, i_2, i_3, j_0 ) \notin \mathcal{J}$ and $( \bm{u}^{(1)},\bm{b}^{(1)} )$, $( \bm{u}^{(2)}, \bm{b}^{(2)} )$ are two solutions of \eqref{eq:general:MHD} on $(0,T)$, and satisfy \eqref{eq:reg} on $(0,T)$. 
	We further assume that 
	\begin{align}
		\Lambda_k^\frac{9}{4}u_l^{(2)}\in L^2(0,T;L^2), \quad \text{if there exist $k,l\in\{1,2,3\}$ s.t. $k\not=l$, $i_k=i_l=k$}, \label{eq:uniqueness:cond:1} \\
		\Lambda_k^\frac{9}{4}b_l^{(2)}\in L^2(0,T;L^2), \quad \text{if there exist $k,l\in\{1,2,3\}$ s.t. $k\not=l$, $i_k=j_0=k$}. \label{eq:uniqueness:cond:2}
	\end{align}
	Then $(\bm{u}^{(1)},\bm{b}^{(1)}) \equiv (\bm{u}^{(2)},\bm{b}^{(2)})$ in $\mathbb{R}^3\times(0,T)$.
\end{theorem}

Our results extend the previous findings of \cite{YJW2019M} in two aspects. Firstly, we allow each component of ${\bm u}$ and ${\bm b}$ to lack dissipation in certain directions. Secondly, we have identified all possible cases where existence and conditional uniqueness can be obtained using this method.
The total number of "good" indices, i.e. the indices which is not in $\mathcal{J}$, is 24. We provide some new cases which ensures the existence of solutions when each component of the velocity fields and magnetic fields lacks dissipation in certain directions. 
However, we did not consider the effects of strong dissipation in a certain direction, such as the vertical direction, as shown in reference \cite{LZ2021} et al.

It is worth noting that, due to the numerous choices for $i_k$, $k=1,2,3$, and $j_0$, the main difficulty of this article lies in the considerable computational workload, as well as providing an analytical proof for any general good indices $i_k$ and $j_0$. The proof of the main theorem relies on estimates of a series of trilinear terms. For instance, when proving the existence of MHD equations, for each chosen $i_k$ and $j_0$, we must estimate $3\times 9$ trilinear terms in equation \eqref{eq:K}, namely $K_1$, $K_2$ and $K_4$, even accounting for symmetry by setting $k=1$ in equation \eqref{eq:K}. When $i_k$ and $j_0$ are arbitrary, it roughly translates to estimating around $3^4 \times 3\times 9$ trilinear terms. Moreover, there are additional estimates required for demonstrating the uniqueness result. In fact, we have categorized estimates of all trilinear terms into four groups, corresponding to Lemma \ref{lem:case:1}-\ref{lem:case:4}. Lemma \ref{lem:case:4} is exclusively utilized for establishing uniqueness. The application of Lemma \ref{lem:case:1}-\ref{lem:case:3} imposes fixed requirements on the indices of the trilinear terms. In the process of deriving the results presented in this article, we initially employed a computer program to identify all good indices for which existence result holds. For each bad index, we outlined those terms for which estimation could not be accomplished using the methods described in this article. Subsequently, we dedicated time to completing an analytical proof for any general good indices $i_k$ and $j_0$.

Define the operator $\widetilde{\Lambda}_{k}^\gamma$ for any $\gamma>0$ by 
\[
\begin{aligned}
	& \widehat {\widetilde{\Lambda}_1^\gamma f( \xi )} := ( \xi_2^2 + \xi_3^2 )^{\gamma/2} \widehat {f( \xi )}, \\
	& \widehat {\widetilde{\Lambda}_2^\gamma f( \xi )} := ( \xi_1^2 + \xi_3^2 )^{\gamma/2} \widehat {f( \xi )}, \\
	& \widehat {\widetilde{\Lambda}_3^\gamma f( \xi )} := ( \xi_1^2 + \xi_2^2 )^{\gamma/2} \widehat {f( \xi )}. 
\end{aligned}
\]
For the following modified Navier-Stokes
\begin{equation}\label{eq:general:NS:2}
	\left\{ 
	\begin{array}{l}
		\partial _t u_1 + \boldsymbol{u} \cdot \nabla u_1 + \nu \widetilde{\Lambda}_{i_1} ^{5/2} u_1 + \partial _1 p = 0 ,\\
		\partial _t u_2 + \boldsymbol{u} \cdot \nabla u_2 + \nu \widetilde{\Lambda}_{i_2} ^{5/2} u_2 + \partial _2 p = 0,\\
		\partial _t u_3 + \boldsymbol{u} \cdot \nabla u_3 + \nu \widetilde{\Lambda}_{i_3} ^{5/2} u_3 + \partial _3 p = 0,\\
		\dive \boldsymbol{u} = 0,\\
		\boldsymbol{u}( x,0 ) = \boldsymbol{u}_0( x ),
	\end{array} 
	\right.
\end{equation}
and modified MHD equations
\begin{equation}\label{eq:general:MHD:2}
	\left\{ 
	\begin{array}{l}
		\partial _t u_1 + \boldsymbol{u} \cdot \nabla u_1 + \nu \widetilde{\Lambda}_{i_1}^{5/2} u_1 + \partial _1 p = \boldsymbol{b} \cdot \nabla b_1, \\
		\partial _t u_2 + \boldsymbol{u} \cdot \nabla u_2 + \nu \widetilde{\Lambda}_{i_2}^{5/2} u_2 + \partial _2 p = \boldsymbol{b} \cdot \nabla b_2, \\
		\partial _t u_3 + \boldsymbol{u} \cdot \nabla u_3 + \nu \widetilde{\Lambda}_{i_3}^{5/2} u_3 + \partial _3 p = \boldsymbol{b} \cdot \nabla b_3, \\
		\pt_t {\bm b} + {\bm u} \cdot \nabla {\bm b} + \eta \widetilde{\Lambda}_{j_0}^{5/2} {\bm b} = {\bm b} \cdot \nabla {\bm u}, \\
		\dive \boldsymbol{u} = 0,\qquad \dive \boldsymbol{b} = 0,\\
		\boldsymbol{u}( x,0 ) = \boldsymbol{u}_0( x ), \qquad \boldsymbol{b}( x,0 ) = \boldsymbol{b}_0( x ),
	\end{array} 
	\right.
\end{equation}
similar results to those in Theorem \ref{thm:existence:NS}-Theorem \ref{thm:uniqueness:MHD} can be readily derived. To be precise, we have the following corollaries. 
\begin{corollary}\label{cor:NS}
	If the indices $(i_1, i_2, i_3) \notin \mathcal{I}$ and $ \boldsymbol{u}_0 \in H^1( \mathbb{R}^3 ) $, then there exists a global solution ${\bm u}$ of \eqref{eq:general:NS:2} satisfying
	\eqref{eq:thm:exis}. Moreover, suppose $\bm{u}^{(1)}$ and $\bm{u}^{(2)}$ are two solutions of \eqref{eq:general:NS:2} satisfying \eqref{eq:thm:exis}, and condition \eqref{eq:thm:uniq} holds for some $T>0$ if there exist $k,l\in\{1,2,3\}$, $k\not=l$ such that $i_k=i_l=k$.
	Then $\bm{u}^{(1)} \equiv \bm{u}^{(2)}$ in $\mathbb{R}^3\times(0,T)$. 
\end{corollary}
\begin{corollary}\label{cor:MHD}
	If $(i_1, i_2, i_3, j_0) \notin \mathcal{J}$ and $(\boldsymbol{u}_0, \boldsymbol{b}_0)\in H^1( \mathbb{R}^3 ) $, then there exists a global solution $({\bm u},{\bm b})$ of \eqref{eq:general:MHD:2} satisfying \eqref{eq:reg}. Moreover, suppose $( \bm{u}^{(1)},\bm{b}^{(1)} )$ and $( \bm{u}^{(2)}, \bm{b}^{(2)} )$ are two solutions of \eqref{eq:general:MHD:2} satisfying \eqref{eq:reg}, and condition \eqref{eq:uniqueness:cond:1} and \eqref{eq:uniqueness:cond:2} hold for some $T>0$. Then $\bm{u}^{(1)} \equiv \bm{u}^{(2)}, \bm{b}^{(1)} \equiv \bm{b}^{(2)}$ in $\mathbb{R}^3\times(0,T)$.
\end{corollary}

Throughout the rest of this paper, we use $\|f\|_{L^2}$ to denote $\|f\|_{{L^2}({\mathbb{R}^3})}$, $\|f\|_{L_{x_i}^2}$ to denote the one-dimensional $ L^2 $-norm (in term of $ x_i $), and $ \| f \|_{L_{x_i x_j}^2} (i \ne j)$ to denote the two-dimensional $ L^2 $-norm (in term of $ x_i $ and $ x_j $). In addition, we use the notion
\[
	\| f \|_{L_{x_i}^r L_{x_j}^qL_{x_k}^p} := \Big\| \big\| \| f \|_{L_{x_k}^p} \big\|_{L_{x_j}^q} \Big\|_{L_{x_i}^r}.
\]

The organization of this paper is as follows. In section 2, we prove the global existence of $ H^1 $-solutions to the Navier-Stokes equations \eqref{eq:general:NS}. We prove in section 3 the uniqueness of $ H^1 $-solutions to \eqref{eq:general:NS} under some extra conditions. Section 4 and 5 are devoted to the proof of existence and uniqueness of the MHD equations \eqref{eq:general:MHD}. The lemmas used in the proof of main results are presented in the appendix.

\section{Global existence and regularity for the Navier-Stokes equations}

{\bf Proof of Theorem \ref{thm:existence:NS}.}
	We first make the zeroth-order energy estimates. Taking the inner product of \eqref{eq:general:NS} with ${\bm u}$, integrating by parts and using $\dive {\bm u} = 0$, we get that
	\begin{align*}
		& \| {\bm u}(t) \|_{L^2}^2 + 2\nu \sum_{j\not=i_1} \int_0^t \| \Lambda_j ^\frac{5}{4} u_1 \| _{L^2}^2 d\tau + 2\nu \sum_{j\not=i_2} \int_0^t \| \Lambda_j^\frac{5}{4} u_2 \| _{L^2}^2d\tau \\ 
		& + 2\nu \sum_{j\not=i_3} \int_0^t \| \Lambda_j^\frac{5}{4} u_3 \| _{L^2}^2d\tau 
		= \| \bm u_0 \|_{L^2}^2.
	\end{align*}
	
	Now we do first-order energy estimates for ${\bm u}$. Applying $\nabla $ to \eqref{eq:general:NS} and then taking the inner product of the resulting equations with $\nabla \bm u$, we obtain that
	\begin{align*}
		& \frac{1}{2} \frac{d}{dt} \| \nabla \bm u \|_{L^2}^2 + \nu \sum_{j \not= i_1} \| \Lambda_j ^\frac{5}{4} \nabla u_1 \|_{L^2}^2 + \nu \sum_{j\not=i_2} \| \Lambda_j^\frac{5}{4} \nabla u_2 \|_{L^2}^2 + \nu \sum_{j\not=i_3} \| \Lambda_j^\frac{5}{4} \nabla u_3 \|_{L^2}^2 \\
		& = -\sum_{i,j=1}^3 \int \partial _i u_j \partial _j u_1 \partial _i u_1 dx - \sum_{i,j=1}^3 \int \partial _i u_j \partial _j u_2 \partial _i u_2 dx 
		- \sum_{i,j=1}^3 \int \partial _i u_j \partial _j u_3 \partial _i u_3 dx \\
		& =:I_1+I_2+I_3.
	\end{align*}
	We only estimate $I_1$ in the following context since $I_2,I_3$ can be argued similarly. Now we separate cases according to the indices $i$ and $j$. Define
	\[
	\begin{array}{l}
		\mathcal{A}_1 := \{(i,j)\in\{1,2,3\} \mid i \not= i_1, j \not=i_1\}, \\
		\mathcal{A}_2 := \{(i,j)\in\{1,2,3\} \mid i \not= i_1, j =i_1\} = \mathcal{A}_{21} \cup \mathcal{A}_{22} \cup \mathcal{A}_{23}, \\
		\mathcal{A}_3 := \{(i,j)\in\{1,2,3\} \mid i = i_1, j \not=i_1\} = \mathcal{A}_{31} \cup \mathcal{A}_{32} \cup \mathcal{A}_{33} \cup \mathcal{A}_{34}, \\
		\mathcal{A}_4 := \{(i,j)\in\{1,2,3\} \mid i = i_1, j =i_1\} = \mathcal{A}_{41} \cup \mathcal{A}_{42}, 
	\end{array}
	\]
	where
	\[
	\begin{array}{ll}
		\mathcal{A}_{21} := \{(i,j) \in \mathcal{A}_2 \mid i_1 = 1 \}, & \\
		\mathcal{A}_{22} := \{(i,j) \in \mathcal{A}_2 \mid i_1 \not= 1, i \not= i_j \}, 
		& \mathcal{A}_{23} := \{(i,j) \in \mathcal{A}_2 \mid i_1 \not= 1, i = i_j \}, \\
		\mathcal{A}_{31} := \{(i,j) \in \mathcal{A}_3 \mid i_1 = 1, i_j \not= 1 \}, 
		& \mathcal{A}_{32} := \{(i,j) \in \mathcal{A}_3 \mid i_1 = 1, i_j = 1 \}, \\
		\mathcal{A}_{33} := \{(i,j) \in \mathcal{A}_3 \mid i_1 \not= 1, j \not= 1 \}, 
		& \mathcal{A}_{34} := \{(i,j) \in \mathcal{A}_3 \mid i_1 \not= 1, j = 1 \}, \\
		\mathcal{A}_{41} := \{(i,j) \in \mathcal{A}_4 \mid i_1 = 1 \}, \quad 
		& \mathcal{A}_{42} := \{(i,j) \in \mathcal{A}_4 \mid i_1 \not= 1 \}. 
	\end{array}
	\]
	
	Below, we will estimate the terms in $I_1$ corresponding to indices in each set $\mathcal{A}_i$ or $\mathcal{A}_{ij}$. Our analysis is mainly categorized into four types: Terms with $(i,j)\in \mathcal{A}_1\cup \mathcal{A}_{21}\cup \mathcal{A}_{22}\cup \mathcal{A}_{31}\cup \mathcal{A}_{33}$ can be directly estimated; Terms with $(i,j)\in \mathcal{A}_{23}$ can be estimated through integration by parts; Terms with $(i,j)\in \mathcal{A}_{32}\cup \mathcal{A}_{41}$ need to be decomposed into two terms by solenoidal condition then estimated; And terms with $(i,j)\in \mathcal{A}_{34}\cup \mathcal{A}_{42}$ need to be combined together by divergence-free condition, then estimated.
	
	\vspace{10pt}
	
	{\myemph The case $(i,j)\in \mathcal{A}_1$. }
	
	Since $i \not=i_1$ and $j \not=i_1$, let $k\in \{1,2,3\}\setminus \{i,i_1\}$, then $\{i,i_1,k\}$ is a permutation of $\{1,2,3\}$. Use $f=\pt_iu_j$, $g=\pt_iu_1$, $h=\pt_ju_1$ in Lemma \ref{lem:case:1}, we have
	\begin{align}
		& |\int \pt_i u_j \pt_ju_1\pt_iu_1dx| \nonumber \\
		& \le C\|\pt_iu_j\|_{L^2} \|\Lambda_i^\frac{1}{4}\pt_i u_1\|_{L^2}^\frac{1}{2} \|\Lambda_i^\frac{1}{4} \Lambda_{i_1}\pt_i u_1\|_{L^2}^\frac{1}{2} \|\Lambda_i^\frac{1}{4}\pt_j u_1\|_{L^2}^\frac{1}{2} \|(\Lambda_i^\frac{5}{4},\Lambda_k^\frac{5}{4})\pt_j u_1\|_{L^2}^\frac{1}{2} \nonumber \\
		& \le \frac{\nu}{100}\|(\Lambda_i^\frac{5}{4},\Lambda_k^\frac{5}{4})\pt_j u_1\|_{L^2}^2 + \frac{\nu}{100}\| \Lambda_i^\frac{5}{4} \pt_{i_1} u_1\|_{L^2}^2+C \|(\Lambda_i^\frac{5}{4},\Lambda_j^\frac{5}{4}) u_1\|^2_{L^2}\|\pt_iu_j\|_{L^2}^2. \label{ineq:A1}
	\end{align}
	The first and second terms on the right hand side of \eqref{ineq:A1} are controllable since $i_1\notin\{i,k\}$. Meanwhile the higher order derivatives $\|(\Lambda_i^\frac{5}{4},\Lambda_j^\frac{5}{4}) u_1\|^2_{L^2}$ in the last term of \eqref{ineq:A1} are integrable since $i_1\notin\{i,j\}$.
	
	{\myemph The case $(i,j)\in \mathcal{A}_{21}$. }
	
	Since $i \not=i_1$ and $i_1=j=1$, let $l\in\{1,2,3\}\setminus \{1,i\}$, then $\{1,i,l\}$ is a permutation of $\{1,2,3\}$. Use $f=\pt_1u_1$, $g=\pt_iu_1$, $h=\pt_iu_1$ in Lemma \ref{lem:case:1}, we have 
	\begin{align}
		& |\int(\pt_iu_1)^2\pt_1u_1dx| \nonumber \\
		& \le C\|\partial_1 u_1\|_{L^2} \|\Lambda_i^\frac{1}{4}\partial_i u_1\|_{L^2}^\frac{1}{2} \|\Lambda_i^\frac{1}{4}\Lambda_1\partial_i u_1\|_{L^2}^\frac{1}{2} \|\Lambda_i^\frac{1}{4}\partial_i u_1\|_{L^2}^\frac{1}{2} \|(\Lambda_i^\frac{5}{4},\Lambda_l^\frac{5}{4})\partial_i u_1\|_{L^2}^\frac{1}{2}  \nonumber \\
		& \le \frac{\nu}{100}\|(\Lambda_i^\frac{5}{4},\Lambda_l^\frac{5}{4})\partial_i u_1\|_{L^2}^2 + \frac{\nu}{100}\|\Lambda_i^\frac{5}{4}\pt_1u_1
		\|_{L^2}^2
		+ C\|\Lambda_i^\frac{5}{4} u_1\|^2_{L^2}\|\partial_1 u_1\|_{L^2}^2. \label{ineq:A21}
	\end{align}
	The first and second terms on the right hand side of \eqref{ineq:A21} are controllable since $i_1=1\notin\{i,l\}$. The higher order derivative $\|\Lambda_i^\frac{5}{4} u_1\|^2_{L^2}$ in the last term of \eqref{ineq:A21} is integrable. 
	
	{\myemph The case $(i,j)\in \mathcal{A}_{22}$. }
	
	Since $i \not=i_1$, $j=i_1\not=1$ and $i\not=i_j$, let $k\in \{1,2,3\}\setminus \{i,i_1\}$, then $\{i,i_1,k\}$ is a permutation of $\{1,2,3\}$. Use $f=\pt_ju_1$, $g=\pt_iu_j$, $h=\pt_iu_1$ in Lemma \ref{lem:case:1}, we have
	\begin{align}
		& \Big| \int\pt_iu_j\pt_ju_1\pt_iu_1dx \Big| \nonumber \\
		& \le C\|\pt_ju_1\|_{L^2} \|\Lambda_i^\frac{1}{4}\pt_i u_j\|_{L^2}^\frac{1}{2} \|\Lambda_i^\frac{1}{4} \Lambda_{i_1}\pt_i u_j\|_{L^2}^\frac{1}{2} \|\Lambda_i^\frac{1}{4}\pt_i u_1\|_{L^2}^\frac{1}{2} \|(\Lambda_i^\frac{5}{4},\Lambda_k^\frac{5}{4})\pt_i u_1\|_{L^2}^\frac{1}{2} \nonumber \\
		& \le \frac{\nu}{100}\|(\Lambda_i^\frac{5}{4},\Lambda_k^\frac{5}{4})\pt_i u_1\|_{L^2}^2 + \frac{\nu}{100}\| \Lambda_i^\frac{5}{4} \pt_{i_1} u_j\|_{L^2}^2+C\|\Lambda_i^\frac{5}{4} u_j\|_{L^2}\|\Lambda_i^\frac{5}{4} u_1\|_{L^2}\|\pt_ju_1\|_{L^2}^2. \label{ineq:A22}
	\end{align}
	The first and second terms on the right hand side of \eqref{ineq:A22} are controllable since $i_1\notin\{i,k\}$ and $ i_j \not= i$. The higher order derivatives in the last term are integrable for the same reason. 
	
	{\myemph The case $(i,j)\in \mathcal{A}_{31}$. }
	
	Since $i=i_1=1$, $j\not=1$ and $i_j\not=1$, let $k\in\{1,2,3\}\setminus \{1,i_j\}$, we have $\{1,i_j,k\}$ is a permutation of $\{1,2,3\}$. Use $f=\pt_1u_1$, $g=\pt_ju_1$, $h=\pt_1u_j$ in Lemma \ref{lem:case:1}, we have
	\begin{align}
		& | \int \pt_1 u_j \pt_j u_1 \pt_1 u_1 dx | \nonumber \\
		& \le C\|\partial_1 u_1\|_{L^2} \|\Lambda_k^\frac{1}{4}\partial_j u_1\|_{L^2}^\frac{1}{2} \|\Lambda_k^\frac{1}{4}\Lambda_{i_j}\partial_j u_1\|_{L^2}^\frac{1}{2} \|\Lambda_k^\frac{1}{4}\partial_1 u_j\|_{L^2}^\frac{1}{2} \|(\Lambda_k^\frac{5}{4},\Lambda_1^\frac{5}{4})\partial_1 u_j\|_{L^2}^\frac{1}{2} \nonumber \\
		& \le \frac{\nu}{100}\|(\Lambda_k^\frac{5}{4},\Lambda_{i_j}^\frac{5}{4})\partial_j u_1\|_{L^2}^2 + \frac{\nu}{100}\|(\Lambda_k^\frac{5}{4},\Lambda_1^\frac{5}{4})\partial_1u_j\|_{L^2}^2 \nonumber \\
		& + C \| ( \Lambda_k^\frac{5}{4}, \Lambda_j^\frac{5}{4}) u_1\|_{L^2} \|(\Lambda_k^\frac{5}{4}, \Lambda_1^\frac{5}{4}) u_j\|_{L^2} \|\partial_1 u_1\|_{L^2}^2. \label{ineq:A31}
	\end{align}
	The first and second terms on the right hand side of \eqref{ineq:A31} are controllable since $i_1=1\notin\{k,i_j\}$ and $i_j\notin\{1,k\}$. Meanwhile the higher order derivatives in the last term are integrable since $i_1\notin\{k,j\}$. 
	
	{\myemph The case $(i,j)\in \mathcal{A}_{33}$.}
	
	Assume that $i=i_j$. Since $j\not=i_1$, $j\not=1$ and $i=i_1\not=1$, we have $i=i_1=i_j\not=1$ and $i, j, 1$ are not equal to each other pairwise. Then $\{i_i,i_j,i_1\}\in \mathcal{I}$, which means this term does not appear in $I_1$. Therefore, it holds that $i\not=i_j$. 
	
	Now $i\not=i_j$, let $k\in\{1,2,3\}\setminus \{i,i_j\}$, we have $\{i,i_j,k\}$ is a permutation of $\{1,2,3\}$. Use $f=\pt_iu_1$, $g=\pt_iu_j$, $h=\pt_ju_1$ in Lemma \ref{lem:case:1}, we have 
	\begin{align}
		& | \int \pt_i u_j \pt_j u_1 \pt_i u_1 dx | \nonumber \\
		& \le C \|\partial_i u_1\|_{L^2} \|\Lambda_k^\frac{1}{4}\partial_i u_j\|_{L^2}^\frac{1}{2} \|\Lambda_k^\frac{1}{4}\Lambda_i\partial_i u_j\|_{L^2}^\frac{1}{2} \|\Lambda_k^\frac{1}{4}\partial_j u_1\|_{L^2}^\frac{1}{2} \|(\Lambda_k^\frac{5}{4},\Lambda_{i_j}^\frac{5}{4})\partial_j u_1\|_{L^2}^\frac{1}{2} \nonumber \\
		& \le \frac{\nu}{100}\|(\Lambda_k^\frac{5}{4}, \Lambda_i^\frac{5}{4}) \partial_i u_j\|_{L^2}^2 + \frac{\nu}{100}\|(\Lambda_k^\frac{5}{4}, \Lambda_{i_j}^\frac{5}{4}) \partial_j u_1\|_{L^2}^2 \nonumber \\
		& + C \| ( \Lambda_k^\frac{5}{4}, \Lambda_i^\frac{5}{4}) u_j\|_{L^2} \|(\Lambda_k^\frac{5}{4}, \Lambda_j^\frac{5}{4}) u_1\|_{L^2} \|\partial_i u_1\|_{L^2}^2. \label{ineq:A33}
	\end{align}
	The first and second terms on the right hand side of \eqref{ineq:A33} are controllable since $i_j\notin\{k,i\}$ and $i_1\notin\{k,i_j\}$. Meanwhile the higher order derivatives in the last term are integrable since $i_1\notin\{k,j\}$ and $i_j\notin\{k,i\}$.
	
	{\myemph The case $(i,j)\in \mathcal{A}_{23}$. }
	
	Since $i=i_j\not=i_1$ and $j=i_1\not=1$, let $k\in \{1,2,3\}\setminus \{i,j\}$, then $\{i,j,k\}$ is a permutation of $\{1,2,3\}$. Integration by parts yields 
	\[
		- \int \pt_iu_j\pt_ju_1\pt_iu_1dx=\int u_1\pt_j\pt_iu_1\pt_iu_jdx+\int u_1\pt_j\pt_iu_j\pt_iu_1dx=:I_a+I_b.
	\]
	
	Use $f=\pt_j\pt_iu_1$, $g=u_1$, $h=\pt_iu_j$ in Lemma \ref{lem:case:3}, we have
	\begin{align}
		|I_a| & \le C\|\Lambda_k^\frac{1}{4}\pt_j\pt_iu_1\|_{L^2} \|\partial_i u_j\|_{L^2}^\frac{2}{5} \|(\Lambda_k^\frac{5}{4},\Lambda_j^\frac{5}{4})\partial_i u_j\|_{L^2}^\frac{3}{5} \|u_1\|_{L^2}^\frac{3}{5} \|\Lambda_i^\frac{5}{4} u_1\|_{L^2}^\frac{2}{5} \nonumber \\
		& \le \frac{\nu}{100}\|(\Lambda_k^\frac{5}{4},\Lambda_j^\frac{5}{4})\partial_i u_j\|_{L^2}^2 + \frac{\nu}{100}\|(\Lambda_k^\frac{5}{4},\Lambda_i^\frac{5}{4})\partial_j u_1\|_{L^2}^2 + C\| u_1\|_{L^2}^3 \|\Lambda_i^\frac{5}{4} u_1\|_{L^2}^2 \|\partial_i u_j\|_{L^2}^2. \label{ineq:A23:1}
	\end{align}
	The first and second terms on the right hand side of \eqref{ineq:A23:1} are controllable since $i_j=i\notin\{k,j\}$ and $i_1\notin\{i,k\}$. Meanwhile the higher order derivative in the last term is integrable. 
	
	Use $f=\pt_j\pt_iu_j$, $g=\pt_iu_1$, $h=u_1$ in	Lemma \ref{lem:case:2}, we have
	\begin{align}
		|I_b| & \le C\|\Lambda_k^\frac{1}{4}\pt_j\pt_iu_j\|_{L^2} \|\Lambda_i^\frac{1}{4}\pt_iu_1\|_{L^2}^\frac{1}{2}\|\Lambda_i^\frac{1}{4}\Lambda_j\pt_iu_1\|_{L^2}^\frac{1}{2} \|u_1\|_{L^2}^\frac{1}{2} \|(\pt_k,\pt_i)u_1\|_{L^2}^\frac{1}{2} \nonumber \\
		& \le \frac{\nu}{100}\|(\Lambda_k^\frac{5}{4},\Lambda_j^\frac{5}{4})\pt_i u_j\|_{L^2}^2 + \frac{\nu}{100}\|\Lambda_i^\frac{5}{4}\pt_ju_1\|_{L^2}^2 + C\| u_1\|_{L^2}^2 \|\Lambda_i^\frac{5}{4} u_1\|_{L^2}^2 \|(\pt_k,\pt_i)u_1\|_{L^2}^2. \label{ineq:A23:2}
	\end{align}
	The first and second terms on the right hand side of \eqref{ineq:A23:2} are controllable since $i_j=i\notin\{k,j\}$ and $i_1\not=i$. Meanwhile the higher order derivative in the last term is integrable. 
	
	{\myemph The case $(i,j)\in \mathcal{A}_{32}$.}
	
	Since $i=i_1=i_j=1\not=j$, let $l\in\{1,2,3\}\setminus \{1,j\}$, we have $\{1,j,l\}$ is a permutation of $\{1,2,3\}$. Due to the solenoidal condition, we have 
	\[
		-\int\pt_1u_j\pt_ju_1\pt_1u_1dx=\int\pt_1u_j\pt_ju_1\pt_ju_jdx+\int\pt_1u_j\pt_ju_1\pt_lu_ldx=:I_a+I_b.
	\]
	
	Use $f=\pt_1u_j$, $g=\pt_ju_1$, $h=\pt_ju_j$ in Lemma \ref{lem:case:1}, we have 
	\begin{align}
		|I_a| & \le C\|\pt_1 u_j\|_{L^2} \|\Lambda_j^\frac{1}{4}\pt_j u_j\|_{L^2}^\frac{1}{2} \|(\Lambda_j^\frac{5}{4},\Lambda_l^\frac{5}{4})\pt_ju_j\|_{L^2}^\frac{1}{2} \|\Lambda_j^\frac{1}{4}\pt_j u_1\|_{L^2}^\frac{1}{2} \|\Lambda_j^\frac{1}{4}\Lambda_1\pt_j u_1\|_{L^2}^\frac{1}{2} \nonumber \\
		& \le \frac{\nu}{100}\|(\Lambda_j^\frac{5}{4},\Lambda_l^\frac{5}{4})\pt_j u_j\|_{L^2}^2 + \frac{\nu}{100}\|\Lambda_j^\frac{5}{4}\pt_1 u_1\|_{L^2}^2 + C\|\Lambda_j^\frac{5}{4} u_1\|_{L^2} \|\Lambda_j^\frac{5}{4} u_j\|_{L^2}\|\pt_1 u_j\|_{L^2}^2. \label{ineq:A32:1}
	\end{align}
	The first and second terms on the right hand side of \eqref{ineq:A32:1} are controllable since $i_j=1\notin\{j,l\}$ and $i_1=1\not=j$. Meanwhile the higher order derivatives in the last term are integrable. 
	
	For $I_b$, we use integration by parts and obtain that 
	\[
		I_b=-\int u_l\pt_l\pt_1u_j\pt_ju_1dx-\int u_l\pt_l\pt_ju_1\pt_1u_jdx=:I_c+I_d.
	\]
	
	Use $f=\pt_l\pt_1u_j$, $g=\pt_ju_1$, $h=u_l$ in Lemma \ref{lem:case:2}, we have
	\begin{align}
		|I_c| & \le C\|\Lambda_l^\frac{1}{4}\pt_l\pt_1u_j\|_{L^2} \|\Lambda_j^\frac{1}{4}\pt_ju_1\|_{L^2}^\frac{1}{2}\|\Lambda_j^\frac{1}{4}\Lambda_1\pt_ju_1\|_{L^2}^\frac{1}{2} \|u_l\|_{L^2}^\frac{1}{2} \|(\pt_l,\pt_j)u_l\|_{L^2}^\frac{1}{2} \nonumber \\
		& \le \frac{\nu}{100}\|\Lambda_l^\frac{5}{4}\partial_1 u_j\|_{L^2}^2 + \frac{\nu}{100}\|\Lambda_j^\frac{5}{4}\pt_1u_1\|_{L^2}^2 + C\| u_l\|_{L^2}^2 \|\Lambda_j^\frac{5}{4} u_1\|_{L^2}^2 \|(\pt_l,\pt_j)u_l\|_{L^2}^2. \label{ineq:A32:2}
	\end{align}
	The first and second terms on the right hand side of \eqref{ineq:A32:2} are controllable since $i_j=i_1=1\notin\{j,l\}$. Meanwhile the higher order derivative in the last term is integrable. 
	
	From $i_1=i_j=1$, we know that $i_l\not=i_j$, since otherwise $i_l=i_j=1$ and $\{i_1,i_l,i_j\}\in \mathcal{I}$. 
	Then Lemma \ref{lem:case:3} is applicable. Use $f=\pt_l\pt_ju_1$, $g=u_l$, $h=\pt_1u_j$ in Lemma \ref{lem:case:3}, we have
	\begin{align}
		|I_d| & \le C\|\Lambda_j^\frac{1}{4}\pt_l\pt_ju_1\|_{L^2} \|\pt_1 u_j\|_{L^2}^\frac{2}{5} \|(\Lambda_j^\frac{5}{4},\Lambda_l^\frac{5}{4})\pt_1 u_j\|_{L^2}^\frac{3}{5} \|u_l\|_{L^2}^\frac{3}{5} \|\Lambda_1^\frac{5}{4} u_l\|_{L^2}^\frac{2}{5} \nonumber \\
		& \le \frac{\nu}{100}\|\Lambda_j^\frac{5}{4}\pt_l u_1\|_{L^2}^2 + \frac{\nu}{100}\|(\Lambda_j^\frac{5}{4},\Lambda_l^\frac{5}{4})\pt_1u_j\|_{L^2}^2 + C\|u_l\|_{L^2}^3 \|\Lambda_1^\frac{5}{4} u_l\|_{L^2}^2 \|\pt_1 u_j\|_{L^2}^2. \label{ineq:A32:3}
	\end{align}
	The first and second terms on the right hand side of \eqref{ineq:A32:3} are controllable since $i_1=1\not=j$ and $i_j=1\notin\{j,l\}$. The higher order derivative in the last term is integrable since $i_l \not= 1$. 
	
	{\myemph The case $(i,j)\in \mathcal{A}_{41}$.}
	
	Now we have $i=j=i_1=1$. Due to the solenoidal condition, we have 
	\begin{align*}
		- \int\pt_iu_j\pt_ju_1\pt_iu_1dx & = -\int (\pt_1u_1)^3dx \\
		& = \int (\pt_1u_1)^2\pt_2u_2dx+\int (\pt_1u_1)^2\pt_3u_3dx=:I_a+I_b.
	\end{align*}
	In order to control $I_a$, we need to discuss based on whether $i_2$ equals one.
	
	(i) $i_2\not=1$. Integration by parts yields that 
	\[
		I_a=-2\int u_2\pt_2\pt_1u_1\pt_1u_1dx.
	\]
	Use $f=\pt_2\pt_1u_1$, $g=u_2$, $h=\pt_1u_1$ in Lemma \ref{lem:case:3}, we have
	\begin{align}
		|I_a| & \le C\|\Lambda_2^\frac{1}{4}\pt_2\pt_1u_1\|_{L^2} \|\partial_1 u_1\|_{L^2}^\frac{2}{5} \|(\Lambda_2^\frac{5}{4},\Lambda_3^\frac{5}{4})\partial_1 u_1\|_{L^2}^\frac{3}{5} \|u_2\|_{L^2}^\frac{3}{5} \|\Lambda_1^\frac{5}{4} u_2\|_{L^2}^\frac{2}{5} \nonumber \\
		& \le \frac{\nu}{100}\|(\Lambda_2^\frac{5}{4},\Lambda_3^\frac{5}{4})\partial_1 u_1\|_{L^2}^2 + \frac{\nu}{100}\|\Lambda_2^\frac{5}{4}\pt_1u_1\|_{L^2}^2 + C\| u_2\|_{L^2}^3 \|\Lambda_1^\frac{5}{4} u_2\|_{L^2}^2 \|\partial_1 u_1\|_{L^2}^2. \label{ineq:A41:1}
	\end{align}
	The first and second terms on the right hand side of \eqref{ineq:A41:1} are controllable since $i_1=1\notin\{2,3\}$. Meanwhile the higher order derivative in the last term is integrable since $i_2\not=1$. 
	
	(ii) $i_2=1$. We obtain since ${\bm u}$ is divergence free that
	\[
		I_a = -\int (\pt_2u_2+\pt_3u_3)\pt_1u_1\pt_2u_2dx =:I_c+I_d.
	\]
	Use $f=\pt_1u_1$, $g=\pt_2u_2$, $h=\pt_2u_2$ in Lemma \ref{lem:case:1}, we have
	\begin{align}
		|I_c| & \le C\|\pt_1u_1\|_{L^2} \|\Lambda_2^\frac{1}{4}\pt_2 u_2\|_{L^2}^\frac{1}{2} \|\Lambda_2^\frac{1}{4} \Lambda_1\pt_2 u_2\|_{L^2}^\frac{1}{2} \|\Lambda_2^\frac{1}{4}\pt_2 u_2\|_{L^2}^\frac{1}{2} \|(\Lambda_2^\frac{5}{4},\Lambda_3^\frac{5}{4})\pt_2 u_2\|_{L^2}^\frac{1}{2} \nonumber \\
		& \le \frac{\nu}{100}\|(\Lambda_2^\frac{5}{4},\Lambda_3^\frac{5}{4})\pt_2 u_2\|_{L^2}^2 + \frac{\nu}{100}\| \Lambda_2^\frac{5}{4} \pt_1 u_2\|_{L^2}^2+C\|\Lambda_2^\frac{5}{4} u_2\|^2_{L^2}\|\pt_1 u_1\|_{L^2}^2. \label{ineq:A41:2}
	\end{align}
	The first and second terms on the right hand side of \eqref{ineq:A41:2} are controllable since $i_2=1\notin\{2,3\}$. And the higher order derivative in the last term is integrable. 
	
	Integration by parts yields 
	\[
		I_d = \int u_1 \pt_1 \pt_2 u_2 \pt_3 u_3 dx + \int u_1 \pt_1 \pt_3 u_3 \pt_2 u_2 dx =: I_e + I_f. 
	\]
	From $i_1 = i_2 = 1$, we know that $i_3 \not= i_1$, since otherwise $i_3 = i_2 = 1$ and $\{ i_1, i_2, i_3 \}\in \mathcal{I}$. Use $f = \pt_1 \pt_2 u_2$, $g = \pt_3 u_3$, $h = u_1$ in Lemma \ref{lem:case:3}, we have
	\begin{align}
		| I_e | & \le C \| \Lambda_2^\frac{1}{4} \pt_1 \pt_2 u_2 \|_{L^2} \| \Lambda_1^\frac{5}{4} \pt_3 u_3 \|_{L^2}^\frac{2}{5} \| \pt_3 u_3 \|_{L^2}^\frac{3}{5} \| u_1 \|_{L^2}^\frac{2}{5} \| ( \Lambda_2^\frac{5}{4}, \Lambda_3^\frac{5}{4} ) u_1 \|_{L^2}^\frac{3}{5} \nonumber \\
		& \le \frac{\nu}{100} \| \Lambda_2^\frac{5}{4} \partial_1 u_2 \|_{L^2}^2 + \frac{\nu}{100} \| \Lambda_1^\frac{5}{4} \partial_3 u_3 \|_{L^2}^2 
		+ C \| ( \Lambda_2^\frac{5}{4}, \Lambda_3^\frac{5}{4} ) u_1 \|_{L^2}^2 \| u_1 \|_{L^2}^\frac{4}{3} \| \partial_3 u_3 \|_{L^2}^2. \label{ineq:A41:3} 
	\end{align}
	The first and second terms on the right hand side of \eqref{ineq:A41:3} are controllable since $i_2=1\not=2$ and $i_3\not=1$. Meanwhile the higher order derivative in the last term is integrable since $i_1 = 1 \notin \{ 2, 3 \}$. 
	
	Use $f = \pt_1 \pt_3 u_3$, $g = \pt_2 u_2$, $h = u_1$ in Lemma \ref{lem:case:2}, we have
	\begin{align}
		| I_f | & \le C \| \Lambda_1^\frac{1}{4} \pt_1 \pt_3 u_3 \|_{L^2} \| \Lambda_2^\frac{1}{4} \pt_2 u_2 \|_{L^2}^\frac{1}{2} \| \Lambda_2^\frac{1}{4} \Lambda_3 \pt_2 u_2 \|_{L^2}^\frac{1}{2} \| u_1 \|_{L^2}^\frac{1}{2} \| ( \pt_1, \pt_2 ) u_1 \|_{L^2}^\frac{1}{2} \nonumber \\
		& \le \frac{\nu}{100} \| \Lambda_1^\frac{5}{4} \partial_3 u_3 \|_{L^2}^2 + \frac{\nu}{100} \| ( \Lambda_2^\frac{5}{4}, \Lambda_3^\frac{5}{4} ) \pt_2 u_2 \|_{L^2}^2 + C \| u_1 \|_{L^2}^2 \| \Lambda_2^\frac{5}{4} u_2 \|_{L^2}^2 \| ( \pt_1,\pt_2 ) u_1 \|_{L^2}^2. \label{ineq:A41:4}
	\end{align}
	The first and second terms on the right hand side of \eqref{ineq:A41:4} are controllable since $i_3\not=1$ and $i_2=1\notin\{2,3\}$. And the higher order derivative in the last term is integrable. 
	
	The estimates of $I_b$ are similar to those of $I_a$, we need to discuss based on whether $i_3$ equals one. If $i_3=1$, we need to use the solenoidal condition again to obtain the boundedness of $I_b$. For the sake of brevity, we omit the details.
	
	{\myemph The case $(i,j)\in \mathcal{A}_{34} \cup \mathcal{A}_{42}$. }
	
	Now we have $i=i_1\not=1$, $j=1$ or $i$. Add the two terms $j=1$ and $j=i$ together to obtain
	\[
		-\int \pt_i u_1\pt_1 u_1\pt_i u_1dx-\int \pt_i u_i\pt_i u_1\pt_i u_1dx=\int (\pt_i u_1)^2\pt_l u_ldx,
	\]
	where $l\in\{1,2,3\}\backslash\{1,i\}$. Integration by parts yields that
	\[
		\int (\pt_i u_1)^2\pt_l u_ldx=-2\int u_l\pt_l\pt_iu_1\pt_iu_1dx.
	\]
	
	Since $i_1=i$, we know that $i_l\not=i$. Use $f=\pt_l\pt_i u_1$, $g=u_l$, $h=\pt_iu_1$ and $\{1,i,l\}$ is a permutation of $\{1,2,3\}$ in Lemma \ref{lem:case:3}, we have
	\begin{align}
		|\int u_l\pt_iu_1\pt_l\pt_iu_1dx|&\le C\|\Lambda_l^\frac{1}{4}\partial_l\pt_i u_1\|_{L^2} \|\partial_i u_1\|_{L^2}^\frac{2}{5} \|(\Lambda_l^\frac{5}{4},\Lambda_1^\frac{5}{4})\partial_i u_1\|_{L^2}^\frac{3}{5} \|u_l\|_{L^2}^\frac{3}{5} \|\Lambda_i^\frac{5}{4} u_l\|_{L^2}^\frac{2}{5} \nonumber \\
		&\le \frac{\nu}{100}\|(\Lambda_l^\frac{5}{4},\Lambda_1^\frac{5}{4}) \pt_i u_1\|_{L^2}^2 + C\|u_l\|_{L^2}^3 \|\Lambda_i^\frac{5}{4} u_l\|_{L^2}^2 \|\partial_i u_1\|_{L^2}^2. \label{ineq:A3442}
	\end{align}
	The first term on the right hand side of \eqref{ineq:A3442} is controllable since $i_1\notin\{1,l\}$. And the higher order derivative in the last term is integrable since $i_l \not= i$. 
	
	\vspace{10pt}
	
	We have completed the estimation of all terms in $I_1$. The other two terms, $I_2$ and $I_3$, can be estimated similarly to $I_1$ and we omit the details for brevity. Incorporating all of these estimates, we obtain that 
	\begin{align*}
		& \frac{d}{dt} \| \nabla {\bm u}\|_{L^2}^2 + \nu \sum_{j\not=i_1} \| \Lambda_j ^\frac{5}{4} \nabla u_1 \| _{L^2}^2 + \nu \sum_{j\not=i_2} \| \Lambda_j^\frac{5}{4}\nabla u_2 \| _{L^2}^2 + \nu \sum_{j\not=i_3} \| \Lambda_j^\frac{5}{4} \nabla u_3 \| _{L^2}^2 \\
		& \le C \Big(\sum_{j\not=i_1} \| \Lambda_j ^\frac{5}{4} u_1 \| _{L^2}^2 + \sum_{j\not=i_2} \| \Lambda_j^\frac{5}{4} u_2 \| _{L^2}^2+ \sum_{j\not=i_3} \| \Lambda_j^\frac{5}{4} u_3 \| _{L^2}^2\Big)(1 + \| {\bm u} \|_{L^2}^3)\|\nabla {\bm u} \|_{L^2}^2.
	\end{align*}
	By Gr\"{o}nwall's inequality, we have that
	\begin{align*}
		& \|\nabla {\bm u}\|_{L^2}^2 + \nu \int_0^t \Big( \sum_{j\not=i_1} \| \Lambda_j^\frac{5}{4} \nabla u_1 \|_{L^2}^2 
		+ \sum_{j\not=i_2} \| \Lambda_j^\frac{5}{4} \nabla u_2 \|_{L^2}^2 + \sum_{j\not=i_3} \| \Lambda_j^\frac{5}{4} \nabla u_3 \|_{L^2}^2 \Big) d \tau\\
		& \le \| \nabla {\bm u}_0\|_{L^2}^2 \exp{\Big( C \int_0^t \big( \sum_{j\not=i_1}\| \Lambda_j^\frac{5}{4} u_1 \|_{L^2}^2 +\sum_{j\not=i_2}\| \Lambda_j^\frac{5}{4} u_2\|_{L^2}^2 + \sum_{j\not=i_3}\| \Lambda_j^\frac{5}{4} u_3\|_{L^2}^2\big)( 1 + \| {\bm u} \|_{L^2}^3)ds \Big) } \\
		& \le \| \nabla {\bm u}_0\|_{L^2}^2 \exp{ \Big( C \big( 1 + \| {\bm u_0} \|_{L^2}^3\big)\int_0^t \big( \sum_{j\not=i_1}\| \Lambda_j^\frac{5}{4} u_1 \|_{L^2}^2 + \sum_{j\not=i_2}\| \Lambda_j^\frac{5}{4} u_2\|_{L^2}^2 + \sum_{j\not=i_3}\| \Lambda_j^\frac{5}{4} u_3\|_{L^2}^2 \big) ds \Big) } \\
		& \le \| \nabla {\bm u}_0\|_{L^2}^2 \exp{ \Big( C \big( 1 + \|{\bm u_0}\|_{L^2}^3\big)\|{\bm u}_0\|_{L^2}^2 \Big) }. 
	\end{align*}
	With this, we have completed the proof for global existence and regularity.
\hfill$\square$

\section{Uniqueness for the Navier-Stokes equations}

{\bf Proof of Theorem \ref{thm:uniqueness:NS}.} 
	It is sufficient to estimate the $L^2$-difference of two solutions,
	\[
		\widetilde{\bm{u}} = \bm{u}^{(1)} - \bm{u}^{(2)},
	\]
	which satisfies
	\begin{equation}\label{eq:NS:different}
		\left\{ 
		\begin{array}{l}
			\partial _t \tilde{u}_1 + {\bm u}^{(1)} \cdot \nabla \tilde{u}_1 + \widetilde{{\bm u}}\cdot \nabla u_1^{(2)}= - \partial _1 \tilde{p} - \nu (\Lambda_1^\frac{5}{2}+\Lambda_2^\frac{5}{2} +\Lambda_3^\frac{5}{2}-\Lambda_{i_1}^\frac{5}{2}) \tilde{u}_1,\\
			\partial _t \tilde{u}_2 + {\bm u}^{(1)} \cdot \nabla \tilde{u}_2 + \widetilde{{\bm u}} \cdot \nabla u_2^{(2)}= - \partial _2 \tilde{p} - \nu (\Lambda_1^\frac{5}{2}+\Lambda_2^\frac{5}{2} +\Lambda_3^\frac{5}{2}-\Lambda_{i_2}^\frac{5}{2}) \tilde{u}_2,\\
			\partial _t \tilde{u}_3 + {\bm u}^{(1)} \cdot \nabla \tilde{u}_3 + \widetilde{{\bm u}} \cdot \nabla u_3^{(2)}= - \partial _3 \tilde{p} - \nu (\Lambda_1^\frac{5}{2}+\Lambda_2^\frac{5}{2} +\Lambda_3^\frac{5}{2}-\Lambda_{i_3}^\frac{5}{2}) \tilde{u}_3,\\
			\dive \widetilde{\bm u} = 0, \\
			\widetilde{\bm u}(x,0) = \widetilde{\bm u}_0 (x).
		\end{array} 
		\right.
	\end{equation}
	Here $ \tilde{p} := p^{(1)} - p^{( 2 )} $ denotes the difference between the corresponding pressures. Taking inner product of \eqref{eq:NS:different} and $\widetilde{\bm u}$, then integrating by parts, we have
	\begin{align*}
		&\hspace{1em}\frac{1}{2}\frac{d}{dt}\| \widetilde{\bm{u}} \|_{L^2}^2 + \nu \sum_{j\not=i_1}\| \Lambda_j ^\frac{5}{4} \tilde{u}_1 \|_{L^2}^2 + \nu \sum_{j\not=i_2}\|\Lambda _j^\frac{5}{4} \tilde{u}_2 \|_{L^2}^2 + \nu\sum_{j\not=i_3}\| \Lambda _j^\frac{5}{4} \tilde{u}_3 \|_{L^2}^2 \\
		&=-\int(\widetilde{{\bm u}}\cdot\nabla){\bm u}^{(2)}\cdot\widetilde{{\bm u}}dx\\
		&=-\int(\widetilde{{\bm u}}\cdot\nabla)u_1^{(2)}\tilde{u}_1dx-\int(\widetilde{{\bm u}}\cdot\nabla)u_2^{(2)}\tilde{u}_2dx-\int(\widetilde{{\bm u}}\cdot\nabla)u_3^{(2)}\tilde{u}_3dx\\
		&=:J_1 + J_2 +J_3.
	\end{align*}
	
	We estimate $J_1$ and write its terms explicitly,
	\[
		J_1 = - \sum_{k=1}^3 \int \tilde{u}_k \pt_k u_1^{(2)} \tilde{u}_1 dx.
	\]
	Now we separate cases according to the indices $k$. Define
	\[
	\begin{array}{l}
		\mathcal{B}_1 := \{k\in\{1,2,3\} \mid i_1 \not= i_k \}, \\
		\mathcal{B}_2 := \{k\in\{1,2,3\} \mid i_1 = i_k \} = \mathcal{B}_{21} \cup \mathcal{B}_{22} \cup \mathcal{B}_{23}, 
	\end{array}
	\]
	where
	\[
	\begin{array}{l}
		\mathcal{B}_{21} := \{k \in \mathcal{B}_2 \mid i_k \not= k \}, \\ 
		\mathcal{B}_{22} := \{k \in \mathcal{B}_2 \mid i_k = k \not=1\},\\ 
		\mathcal{B}_{23} := \{k \in \mathcal{B}_2 \mid i_k = k =1\}.
	\end{array}
	\]
	
	\setcounter{mystep}{0}
	{\myemph The case $k\in \mathcal{B}_1$. }
	
	Since $i_1\not=i_k$, let $l\in\{1,2,3\}\setminus \{i_1,i_k\}$, we have $\{i_1,i_k,l\}$ is a permutation of $\{1,2,3\}$. 
	
	Use $f=\tilde{u}_1$, $g=\partial_k u_1^{(2)}$, $h=\tilde{u}_k$ in Lemma \ref{lem:case:4}, we have
	\begin{align}
		& |\int \tilde{u}_k\pt_ku_1^{(2)}\tilde{u}_1dx| \nonumber \\
		& \le C\|\tilde{u}_1\|_{L^2}^\frac{3}{5} \|\Lambda_l^\frac{5}{4}\tilde{u}_1\|_{L^2}^\frac{2}{5} \|\partial_k u_1^{(2)}\|_{L^2}^\frac{3}{5} \|\Lambda_{i_k}^\frac{5}{4}\partial_k u_1^{(2)}\|_{L^2}^\frac{2}{5} \|\tilde{u}_k\|_{L^2}^\frac{3}{5} \|\Lambda_{i_1}^\frac{5}{4}\tilde{u}_k\|_{L^2}^\frac{2}{5} \nonumber \\
		& \le \frac{\nu}{100}\|\Lambda_l^\frac{5}{4}\tilde{u}_1\|_{L^2}^2 + \frac{\nu}{100}\|\Lambda_{i_1}^\frac{5}{4}\tilde{u}_k\|_{L^2}^2 + C\|\tilde{u}_1\|_{L^2} \|\tilde{u}_k\|_{L^2} \|\partial_k u_1^{(2)}\|_{L^2} \|\Lambda_{i_k}^\frac{5}{4}\partial_k u_1^{(2)}\|_{L^2}^\frac{2}{3}. \label{ineq:B1}
	\end{align}
	The first and second terms on the right hand side of \eqref{ineq:B1} are controllable since $i_1\not=l$ and $i_k\not=i_1$. From the conditions of Theorem \ref{thm:uniqueness:NS}, we know that the first order derivative in the last term is always bounded, and the higher order derivative $\|\Lambda_{i_k}^\frac{5}{4}\partial_k u_1^{(2)}\|_{L^2}$ is integrable in any finite time interval since $i_1\not=i_k$.   
	
	{\myemph The case $k\in \mathcal{B}_{21}$. }
	
	Since $i_1=i_k\not=k$, let $l\in\{1,2,3\}\setminus \{k,i_1\}$, we have $\{k,l,i_1\}$ is a permutation of $\{1,2,3\}$. Use $f=\tilde{u}_k$, $g=\pt_ku_1^{(2)}$, $h=\tilde{u}_1$ in Lemma \ref{lem:case:1}, we have 
	\[
		|\int \tilde{u}_k\pt_ku_1^{(2)}\tilde{u}_1dx| \le C \| \tilde{u}_k \|_{L^2} \|\Lambda_k^\frac{1}{4}\pt_ku_1^{(2)}\|_{L^2}^\frac{1}{2}\|\Lambda_k^\frac{1}{4}\Lambda_{i_1}\pt_ku_1^{(2)}\|_{L^2}^\frac{1}{2}\|\Lambda_k^\frac{1}{4}\tilde{u}_1\|_{L^2}^\frac{1}{2}\|(\Lambda_k^\frac{5}{4},\Lambda_l^\frac{5}{4})\tilde{u}_1\|_{L^2}^\frac{1}{2}.
	\]
	Using the interpolation inequality 
	\[
		\|\Lambda_k^\frac{1}{4}\tilde{u}_1\|_{L_{x_k}^2}\le C\|\tilde{u}_1\|_{L_{x_k}^2}^\frac{4}{5}\|\Lambda_k^\frac{5}{4}\tilde{u}_1\|_{L_{x_k}^2}^\frac{1}{5},
	\]
	we obtain that 
	\begin{align}
		& |\int \tilde{u}_k\pt_ku_1^{(2)}\tilde{u}_1dx|
		\le C \| \tilde{u}_k \|_{L^2} \| \tilde{u}_1 \|_{L^2}^\frac{2}{5} \| ( \Lambda_k^\frac{5}{4}, \Lambda_l^\frac{5}{4} ) \tilde{u}_1 \|_{L^2}^\frac{3}{5} \| \Lambda_k^\frac{5}{4} u_1^{(2)}\|_{L^2}^\frac{1}{2}\|\Lambda_k^\frac{5}{4}\pt_{i_1}u_1^{(2)}\|_{L^2}^\frac{1}{2} \nonumber \\
		& \le \frac{\nu}{100} \|(\Lambda_k^\frac{5}{4},\Lambda_l^\frac{5}{4})\tilde{u}_1\|_{L^2}^2+C\|\tilde{u}_k\|_{L^2}^\frac{10}{7}\|\tilde{u}_1\|_{L^2}^\frac{4}{7}\|\Lambda_k^\frac{5}{4}u_1^{(2)}\|_{L^2}^\frac{5}{7}\|\Lambda_k^\frac{5}{4}\pt_{i_1}u_1^{(2)}\|_{L^2}^\frac{5}{7}. \label{ineq:B21}
	\end{align}
	The first term on the right hand side of \eqref{ineq:B21} is controllable since $i_1\notin\{k,l\}$. The last term is integrable in any finite time interval since $i_1 \not= k$.  
	
	{\myemph The case $k\in \mathcal{B}_{22}$. }
	
	Now we have $i_1=i_k=k\not=1$. Use $f=\tilde{u}_k$, $g=\partial_k u_1^{(2)}$, $h=\tilde{u}_1$ and $\{i,j,k\}$ is a permutation of $\{1,2,3\}$ in Lemma \ref{lem:case:4}, we have
	\begin{align}
		& |\int \tilde{u}_k \partial_k u_1^{(2)} \tilde{u}_1 dx| 
		\le C\|\tilde{u}_k\|_{L^2}^\frac{3}{5} \|\Lambda_{j}^\frac{5}{4}\tilde{u}_k\|_{L^2}^\frac{2}{5} \|\partial_k u_1^{(2)}\|_{L^2}^\frac{3}{5} \|\Lambda_k^\frac{5}{4}\partial_k u_1^{(2)}\|_{L^2}^\frac{2}{5} \|\tilde{u}_1\|_{L^2}^\frac{3}{5} \|\Lambda_i^\frac{5}{4}\tilde{u}_1\|_{L^2}^\frac{2}{5} \nonumber \\
		& \le \frac{\nu}{100}\|\Lambda_i^\frac{5}{4}\tilde{u}_1\|_{L^2}^2 + \frac{\nu}{100}\|\Lambda_j^\frac{5}{4}\tilde{u}_k\|_{L^2}^2 + C\|\tilde{u}_k\|_{L^2} \|\tilde{u}_1\|_{L^2} \|\partial_k u_1^{(2)}\|_{L^2} \|\Lambda_k^\frac{5}{4}\partial_k u_1^{(2)}\|_{L^2}^\frac{2}{3}. \label{ineq:B22}
	\end{align}
	The first and second terms on the right hand side of \eqref{ineq:B22} are controllable since $i_1\not=i$ and $i_k\not=j$. The last term $\|\Lambda_k^\frac{5}{4}\partial_k u_1^{(2)}\|$ is bounded according to the assumption in Theorem \ref{thm:uniqueness:NS}. 
	
	{\myemph The case $k\in \mathcal{B}_{23}$.}
	
	Now we have $i_1=i_k=k=1$. Due to the solenoidal condition, we have 
	\[
		-\int\tilde{u}_1\pt_1u_1^{(2)}\tilde{u}_1dx=\int\tilde{u}_1\pt_2u_2^{(2)}\tilde{u}_1dx+\int\tilde{u}_1\pt_3u_3^{(2)}\tilde{u}_1dx=:I_a+I_b.
	\]
	
	Suppose $i_2\not=1$. Let $l\in\{1,2,3\}\setminus \{i_1,i_2\}$, we have $\{i_1,i_2,l\}$ is a permutation of $\{1,2,3\}$. Use $f=\tilde{u}_1$, $g=\partial_2 u_2^{(2)}$, $h=\tilde{u}_1$ in Lemma \ref{lem:case:4}, we have
	\begin{align}
		|I_a|
		& \le C\|\tilde{u}_1\|_{L^2}^\frac{3}{5} \|\Lambda_{i_2}^\frac{5}{4}\tilde{u}_1\|_{L^2}^\frac{2}{5} \|\partial_2 u_2^{(2)}\|_{L^2}^\frac{3}{5} \|\Lambda_{i_1}^\frac{5}{4}\partial_2 u_2^{(2)}\|_{L^2}^\frac{2}{5} \|\tilde{u}_1\|_{L^2}^\frac{3}{5} \|\Lambda_l^\frac{5}{4}\tilde{u}_1\|_{L^2}^\frac{2}{5} \nonumber \\
		& \le \frac{\nu}{100}\|\Lambda_{i_2}^\frac{5}{4}\tilde{u}_1\|_{L^2}^2 + \frac{\nu}{100}\|\Lambda_l^\frac{5}{4}\tilde{u}_1\|_{L^2}^2 + C\|\tilde{u}_1\|_{L^2}^2 \|\partial_2 u_2^{(2)}\|_{L^2} \|\Lambda_{i_1}^\frac{5}{4}\partial_2 u_2^{(2)}\|_{L^2}^\frac{2}{3}. \label{ineq:B23:1}
	\end{align}
	The first and second terms on the right hand side of \eqref{ineq:B23:1} are controllable since $i_1\notin\{l,i_2\}$. The first order derivative in the last term is bounded, and the higher order derivative $\|\Lambda_{i_1}^\frac{5}{4}\partial_2 u_2^{(2)}\|_{L^2}$ is integrable in any finite time interval since $i_2\not=i_1$.   
	
	Suppose $i_2=1$. Use $f=\tilde{u}_1$, $g=\pt_2u_2^{(2)}$, $h=\tilde{u}_1$ in Lemma \ref{lem:case:1}, we have 
	\[
		|I_a| \le C \| \tilde{u}_1 \|_{L^2} \| \Lambda_2^\frac{1}{4} \pt_2 u_2^{(2)} \|_{L^2}^\frac{1}{2} \| \Lambda_2^\frac{1}{4} \Lambda_1 \pt_2 u_2^{(2)} \|_{L^2}^\frac{1}{2} \| \Lambda_2^\frac{1}{4} \tilde{u}_1 \|_{L^2}^\frac{1}{2} \| ( \Lambda_2^\frac{5}{4}, \Lambda_3^\frac{5}{4} ) \tilde{u}_1 \|_{L^2}^\frac{1}{2}.
	\]
	By inserting the following interpolation inequality into the above expression
	\[
		\| \Lambda_2^\frac{1}{4} \tilde{u}_1 \|_{L_{x_2}^2} \le C \| \tilde{u}_1 \|_{L_{x_2}^2}^\frac{4}{5} \| \Lambda_2^\frac{5}{4} \tilde{u}_1 \|_{L_{x_2}^2}^\frac{1}{5},
	\]
	we obtain that
	\begin{align}
		| I_a | & \le C \| \tilde{u}_1 \|_{L^2}^\frac{7}{5} \| ( \Lambda_2^\frac{5}{4}, \Lambda_3^\frac{5}{4} ) \tilde{u}_1 \|_{L^2}^\frac{3}{5} \| \Lambda_2^\frac{5}{4} u_2^{(2)} \|_{L^2}^\frac{1}{2} \| \Lambda_2^\frac{5}{4} \pt_1 u_2^{(2)} \|_{L^2}^\frac{1}{2} \nonumber \\
		& \le \frac{\nu}{100} \| ( \Lambda_2^\frac{5}{4}, \Lambda_3^\frac{5}{4} ) \tilde{u}_1 \|_{L^2}^2 + C \| \tilde{u}_1 \|_{L^2}^2 \| \Lambda_2^\frac{5}{4} u_2^{(2)} \|_{L^2}^\frac{5}{7} \| \Lambda_2^\frac{5}{4} \pt_1 u_2^{(2)} \|_{L^2}^\frac{5}{7}. \label{ineq:B23:2}
	\end{align}
	The first term on the right hand side of \eqref{ineq:B23:2} is controllable since $i_1\notin\{2,3\}$. The last term is integrable in any finite time interval since $i_2 \not= 2$. 
	The estimates of $I_b$ are similar to those of $I_a$, we need to discuss based on whether $i_3$ equals one. The details are omitted for the sake of brevity. 
	
	\vspace{10pt}
	
	We have finished estimating all terms in $J_1$. The other two terms $J_2$ and $J_3$ can be similarly estimated as $J_1$ and we omit the details. Putting all these estimates together, we have that
	\[
		\frac{d}{dt} \| \widetilde{\bm{u}} \|_{L^2}^2 + \nu \sum_{j\not=i_1} \| \Lambda_j ^\frac{5}{4} \tilde{u}_1 \|_{L^2}^2 + \nu \sum_{j\not=i_2} \| \Lambda _j^\frac{5}{4} \tilde{u}_2 \|_{L^2}^2 + \nu \sum_{j\not=i_3} \| \Lambda _j^\frac{5}{4} \tilde{u}_3 \|_{L^2}^2 \le A(t) \| \widetilde{\bm{u}} \|_{L^2}^2,
	\]
	where 
	\begin{align*}
		A(t) & = C \big( \sum_{j\not=i_1} \| \Lambda_j ^\frac{5}{4} u_1^{(2)} \|_{L^2}^\frac{5}{7} \| \Lambda_j ^\frac{5}{4} \nabla u_1^{(2)} \|_{L^2}^\frac{5}{7} 
		+ \sum_{j\not=i_2}\| \Lambda_j ^\frac{5}{4} u_2^{(2)} \|_{L^2}^\frac{5}{7}\| \Lambda_j ^\frac{5}{4} \nabla u_2^{(2)} \|_{L^2}^\frac{5}{7} \\
		& + \sum_{j\not=i_3} \| \Lambda_j ^\frac{5}{4} u_3^{(2)} \|_{L^2}^\frac{5}{7} \| \Lambda_j ^\frac{5}{4} \nabla u_3^{(2)}\|_{L^2}^\frac{5}{7}
		+ \sum_{j\not=i_1} \| \nabla u_1^{(2)} \|_{L^2} \| \Lambda_j ^\frac{5}{4} \nabla u_1^{(2)} \|_{L^2}^\frac{2}{3} \\
		& + \sum_{j\not=i_2} \| \nabla u_2^{(2)} \|_{L^2} \| \Lambda_j ^\frac{5}{4} \nabla u_2^{(2)} \|_{L^2}^\frac{2}{3}
		+ \sum_{j\not=i_3} \| \nabla u_3^{(2)} \|_{L^2} \| \Lambda_j ^\frac{5}{4} \nabla u_3^{(2)}\|_{L^2}^\frac{2}{3}\\
		& + \sum_{k\not=1} \delta_{ki_1}\delta_{ki_k} \| \nabla u_1^{(2)} \|_{L^2} \|\Lambda_k^\frac{9}{4} u_1^{(2)} \|_{L^2}^\frac{2}{3} + \sum_{k\not=2} \delta_{ki_2} \delta_{ki_k} \| \nabla u_2^{(2)} \|_{L^2} \| \Lambda_k^\frac{9}{4} u_2^{(2)} \|_{L^2}^\frac{2}{3}\\
		& + \sum_{k\not=3} \delta_{ki_3} \delta_{ki_k} \| \nabla u_3^{(2)} \|_{L^2} \| \Lambda_k^\frac{9}{4} u_3^{(2)} \|_{L^2}^\frac{2}{3} 
		\big),
	\end{align*}
	and	$\delta_{ij}$ denotes the Kronecker delta. 
	The regularity assumptions on ${\bm u}^{(2)}$ in \eqref{eq:thm:exis} and \eqref{eq:thm:uniq} ensure the integrability of $A(t)$ on $(0,T)$ for any fixed $T>0$.
	The Gr\"{o}nwall's inequality then implies that $\bm{u}^{(1)} \equiv \bm{u}^{(2)}$ if $\bm{u}_0^{(1)} \equiv \bm{u}_0^{(2)}$. This completes the proof of uniqueness.
\hfill$\square$

Corollary \ref{cor:NS} can be readily derived through Theorem \ref{thm:existence:NS} and Theorem \ref{thm:uniqueness:NS}.

\section{Global existence and regularity for the MHD equations}

{\bf Proof of Theorem \ref{thm:existence:MHD}.} 
	Taking the inner product of \eqref{eq:general:MHD} with $({\bm u}, {\bm b})$, integrating by parts and using $\dive {\bm u} = \dive {\bm b} = 0$, we get
	\begin{align*}
		& \| ( {\bm u}, {\bm b})(t) \|_{L^2}^2 
		+ 2 \nu \sum_{j\not=i_1} \int_0^t \| \Lambda_j ^\frac{5}{4} u_1 \| _{L^2}^2 d\tau 
		+ 2 \nu \sum_{j\not=i_2} \int_0^t \| \Lambda_j^\frac{5}{4} u_2 \| _{L^2}^2d\tau \\
		& + 2\nu \sum_{j\not=i_3} \int_0^t \| \Lambda_j^\frac{5}{4} u_3 \| _{L^2}^2d\tau + 2\eta \sum_{k=1}^3 \sum_{j \not= j_0} \int_0^t \| \Lambda_j ^\frac{5}{4} b_k \| _{L^2}^2 d\tau 
		= \| ({\bm u}_0, {\bm b}_0) \|_{L^2}^2. 
	\end{align*}
	Applying first order derivative to \eqref{eq:general:MHD} and then taking the inner product of the resulting equations with $(\nabla {\bm u},\nabla {\bm b})$, we obtain that 
	\begin{equation}\label{eq:K}
		\begin{aligned}
		&\; \frac{1}{2}\frac{d}{dt}\| (\nabla {\bm u}, \nabla {\bm b}) \|_{L^2}^2+ \nu \sum_{j\not=i_1}\| \Lambda_j ^\frac{5}{4} \nabla u_1 \|_{L^2}^2 + \nu \sum_{j\not=i_2}\| \Lambda_j ^\frac{5}{4} \nabla u_2 \|_{L^2}^2 \\
		& + \nu \sum_{j\not=i_3}\| \Lambda_j ^\frac{5}{4} \nabla u_3 \|_{L^2}^2 
		+ \eta \sum_{k=1}^3 \sum_{j \not= j_0}\| \Lambda_j ^\frac{5}{4} \nabla b_k \|_{L^2}^2 \\
		& = \int \big[ - \nabla ( {\bm u} \cdot \nabla {\bm u} ) \cdot \nabla {\bm u}
		+ \nabla ( {\bm b} \cdot \nabla {\bm b} ) \cdot \nabla {\bm u} 
		+ \nabla ( {\bm b} \cdot \nabla {\bm u} ) \cdot \nabla {\bm b} 
		- \nabla ( {\bm u} \cdot \nabla {\bm b} ) \cdot \nabla {\bm b} \big] dx  \\
		& = \int \Big[ - \partial _i u_j \partial _j u_k \partial _i u_k
		+ \partial_i b_j \partial_j b_k \partial_i u_k 
		+ \partial_i b_j \partial_j u_k \partial_i b_k 
		- \partial _i u_j \partial _j b_k \partial _i b_k \Big] dx  \\
		& =: K_1 + K_2 + K_3 + K_4.
		\end{aligned}
	\end{equation}
	
	The estimates for $K_1$ have been obtained in section 3 for the Navier-Stokes equations. As $K_2$ and $K_3$ exhibit symmetry with respect to $i$ and $j$, we will present estimations for $K_2$ and $K_4$ in the next two subsections, respectively. 
	
	\subsection{The estimates of $K_2$}
	
	The items in $K_2$ can be expanded as
	\begin{align*}
		K_2 = & \int \partial _i b_j \partial _j b_k \partial _i u_k dx
		= \sum\limits_{i,j = 1}^{\rm 3} \int \partial _i b_j\partial _j b_1\partial _i u_1 dx 
		+ \sum\limits_{i,j = 1}^{\rm 3} \int \partial _i b_j\partial _j b_2\partial _i u_2 dx \\
		& + \sum\limits_{i,j = 1}^{\rm 3} \int \partial _i b_j\partial _j b_3\partial _i u_3 dx=: K_{21} + K_{22} + K_{23}.
	\end{align*}
	We only estimate $K_{21}$ in the following context since $K_{22}$, $K_{23}$ can be argued similarly. Now we separate cases according to the indices $i$ and $j$. Define
	\begin{align*}
		\mathcal{D}_1 := & \{(i,j)\in\{1,2,3\} \mid i \not= i_1, j  \not= j_0 \},\\
		\mathcal{D}_2 := & \{(i,j)\in\{1,2,3\} \mid i \not= i_1, j = j_0 \} = \mathcal{D}_{21} \cup \mathcal{D}_{22} \cup \mathcal{D}_{23} \cup \mathcal{D}_{24}, \\
		\mathcal{D}_3 := & \{(i,j)\in\{1,2,3\} \mid i = i_1, j \not= j_0\} = \mathcal{D}_{31} \cup \mathcal{D}_{32}, \\
		\mathcal{D}_4 := & \{(i,j)\in\{1,2,3\} \mid i = i_1, j = j_0\} = \mathcal{D}_{41} \cup \mathcal{D}_{42}, 
	\end{align*}
	where
	\[
	\begin{array}{l}
		\mathcal{D}_{21} := \{(i,j) \in \mathcal{D}_2 \mid j_0 \not= i\}, \\
		\mathcal{D}_{22} := \{(i,j) \in \mathcal{D}_2 \mid j_0 = i, j_0=1 \}, \\
		\mathcal{D}_{23} := \{(i,j) \in \mathcal{D}_2 \mid j_0 = i, j_0\not=1, i_1 = 1 \}, \\
		\mathcal{D}_{24} := \{(i,j) \in \mathcal{D}_2 \mid j_0 = i, j_0\not=1, i_1 \not= 1 \}, 
	\end{array}
	\]
	\[
	\begin{array}{l}
		\mathcal{D}_{31} := \{(i,j) \in \mathcal{D}_3 \mid j_0 \not= i\}, \\
		\mathcal{D}_{32} := \{(i,j) \in \mathcal{D}_3 \mid j_0 = i \}, \\
	\end{array}
	\]
	\[
	\begin{array}{l}
		\mathcal{D}_{41} := \{(i,j) \in \mathcal{D}_4 \mid j_0 = i_1 \}, \\
		\mathcal{D}_{42} := \{(i,j) \in \mathcal{D}_4 \mid j_0 \not= i_1 \}. 
	\end{array}
	\]
	
	\setcounter{mystep}{0}
	{\myemph The case $(i,j)\in \mathcal{D}_1$. }
	
	Now we have $i_1 \not= i$ and $j_0 \not= j$. Assume $i_1=j_0$. According to the definition of set $\mathcal{J}$, now we have $i_1 = j_0 = 1$. Since $i_1 = 1 \not= i$, let $l \in \{1,2,3\}\setminus \{ 1, i \}$, then $\{ 1, i, l \}$ is a permutation of $\{ 1, 2, 3 \}$. Use $ f = \pt_i b_j $, $ g = \pt_j b_1 $, $ h = \pt_i u_1 $ in Lemma \ref{lem:case:1}, we have 
	\begin{align}
		& | \int \pt_i b_j \pt_j b_1 \pt_i u_1 dx | \nonumber \\
		& \le C \| \partial_i b_j \|_{L^2} \| \Lambda_i^\frac{1}{4}\partial_i u_1 \|_{L^2}^\frac{1}{2} \| \Lambda_i^\frac{1}{4} \Lambda_1 \partial_i u_1 \|_{L^2}^\frac{1}{2} \| \Lambda_i^\frac{1}{4} \partial_j b_1 \|_{L^2}^\frac{1}{2} \|( \Lambda_i^\frac{5}{4}, \Lambda_l^\frac{5}{4} ) \partial_j b_1 \|_{L^2}^\frac{1}{2} \nonumber \\
		& \le \frac{\eta}{100} \| ( \Lambda_i^\frac{5}{4}, \Lambda_l^\frac{5}{4} ) \partial_j b_1 \|_{L^2}^2 + \frac{\nu}{100} \| \Lambda_i^\frac{5}{4} \pt_1 u_1 \|_{L^2}^2 + C \| \Lambda_i^\frac{5}{4} u_1 \|_{L^2} \| ( \Lambda_i^\frac{5}{4}, \Lambda_j^\frac{5}{4} ) b_1 \|_{L^2} \| \partial_i b_j \|_{L^2}^2. \label{ineq:D1:1}
	\end{align}
	The first and second terms on the right hand side of \eqref{ineq:D1:1} are controllable since $j_0 = 1 \notin \{ i, l \}$ and $i_1 = 1 \not= i$. Meanwhile the higher order derivatives in the last term are integrable since $i_1 \not= i$ and $j_0 \notin \{ i, j \}$. 
	
	Assume $ i_1 \not= j_0 $, let $ l \in\{ 1, 2, 3 \} \setminus \{ i_1, j_0 \}$, then $\{ i_1, j_0, l \}$ is a permutation of $ \{ 1, 2, 3 \}$. Use $ f = \pt_i b_j $, $ g = \pt_j b_1 $, $ h = \pt_i u_1 $ in Lemma \ref{lem:case:1}, we have 
	\begin{align}
		& | \int \pt_i b_j \pt_j b_1 \pt_i u_1 dx | \nonumber \\
		& \le C \| \pt_i b_j \|_{L^2} \| \Lambda_l^\frac{1}{4} \pt_i u_1 \|_{L^2}^\frac{1}{2} \| \Lambda_l^\frac{1}{4} \Lambda_{j_0} \partial_i u_1 \|_{L^2}^\frac{1}{2} \| \Lambda_l^\frac{1}{4} \partial_j b_1 \|_{L^2}^\frac{1}{2} \| ( \Lambda_l^\frac{5}{4}, \Lambda_{i_1}^\frac{5}{4} ) \partial_j b_1 \|_{L^2}^\frac{1}{2} \nonumber \\
		& \le \frac{\eta}{100} \| ( \Lambda_l^\frac{5}{4}, \Lambda_{i_1}^\frac{5}{4} ) \partial_j b_1 \|_{L^2}^2 + \frac{\nu}{100} \| ( \Lambda_l^\frac{5}{4}, \Lambda_{j_0}^\frac{5}{4} ) \partial_i u_1 \|_{L^2}^2 \nonumber \\
		& + C \| ( \Lambda_l^\frac{5}{4}, \Lambda_i^\frac{5}{4}) u_1 \|_{L^2} \| ( \Lambda_l^\frac{5}{4}, \Lambda_j^\frac{5}{4} ) b_1 \|_{L^2} \| \partial_i b_j \|_{L^2}^2. \label{ineq:D1:2}
	\end{align}
	The first and second terms on the right hand side of \eqref{ineq:D1:2} are controllable since $j_0 \notin \{ i_1, l \} $ and $i_1 \notin\{ j_0, l \}$. Meanwhile the higher order derivatives in the last term are integrable since $i_1\notin\{l,i\}$ and $j_0 \notin \{ l, j \}$.
	
	{\myemph The case $(i,j)\in \mathcal{D}_{21}$. }
	
	Now we have $i_1 \not= i$, $j_0 = j \not= i$. Assume $i_1 \not= j_0$, let $l \in \{ 1, 2, 3 \} \setminus \{ i_1, j_0 \}$, then $\{ i_1, j_0, l \}$ is a permutation of $\{ 1, 2, 3 \}$. Use $f = \pt_j b_1$, $g = \pt_i b_j$, $h = \pt_i u_1$ in Lemma \ref{lem:case:1}, we have 
	\begin{align}
		& | \int \pt_i b_j \pt_j b_1 \pt_i u_1 dx | \nonumber \\
		& \le C \| \partial_j b_1 \|_{L^2} \| \Lambda_l^\frac{1}{4} \partial_i b_j \|_{L^2}^\frac{1}{2} \| \Lambda_l^\frac{1}{4} \Lambda_{i_1} \partial_i b_j \|_{L^2}^\frac{1}{2} \|\Lambda_l^\frac{1}{4} \partial_i u_1 \|_{L^2}^\frac{1}{2} \| ( \Lambda_l^\frac{5}{4}, \Lambda_{j_0}^\frac{5}{4} ) \partial_i u_1 \|_{L^2}^\frac{1}{2} \nonumber \\
		& \le \frac{\eta}{100} \| ( \Lambda_l^\frac{5}{4}, \Lambda_{i_1}^\frac{5}{4} ) \partial_i b_j \|_{L^2}^2 + \frac{\nu}{100} \| ( \Lambda_l^\frac{5}{4}, \Lambda_{j_0}^\frac{5}{4} ) \partial_i u_1 \|_{L^2}^2 \nonumber \\
		& + C \| ( \Lambda_l^\frac{5}{4}, \Lambda_i^\frac{5}{4} ) b_j \|_{L^2} \| ( \Lambda_l^\frac{5}{4}, \Lambda_i^\frac{5}{4} ) u_1 \|_{L^2} \| \partial_j b_1 \|_{L^2}^2. \label{ineq:D21:1}
	\end{align}
	The first and second terms on the right hand side of \eqref{ineq:D21:1} are controllable since $j_0 \notin \{l,i_1\}$ and $i_1\notin\{j_0, l\}$. Meanwhile the higher order derivatives in the last term are integrable since $j_0 \notin\{l,i\}$ and $i_1\notin\{l,i\}$. 
	
	Assume $i_1 = j_0$, since $i_1\not=i$, let $l\in\{1,2,3\}\setminus \{i,i_1\}$, then $\{i,i_1,l\}$ is a permutation of $\{1,2,3\}$. Use $f=\pt_jb_1$, $g=\pt_ib_j$, $h=\pt_iu_1$ in Lemma \ref{lem:case:1}, we have 
	\begin{align}
		& |\int \pt_i b_j \pt_j b_1 \pt_i u_1 dx | \nonumber \\
		& \le C\|\partial_j b_1\|_{L^2} \|\Lambda_i^\frac{1}{4}\partial_i b_j\|_{L^2}^\frac{1}{2} \|\Lambda_i^\frac{1}{4}\Lambda_{i_1}\partial_i b_j\|_{L^2}^\frac{1}{2} \|\Lambda_i^\frac{1}{4}\partial_i u_1\|_{L^2}^\frac{1}{2} \|(\Lambda_i^\frac{5}{4},\Lambda_l^\frac{5}{4})\partial_i u_1\|_{L^2}^\frac{1}{2} \nonumber \\
		& \le \frac{\eta}{100}\|\Lambda_i^\frac{5}{4}\partial_{i_1} b_j\|_{L^2}^2 + \frac{\nu}{100}\|(\Lambda_i^\frac{5}{4},\Lambda_l^\frac{5}{4})\partial_i u_1
		\|_{L^2}^2+ C\|\Lambda_i^\frac{5}{4} b_j\|_{L^2}\|\Lambda_i^\frac{5}{4} u_1\|_{L^2}\|\partial_j b_1\|_{L^2}^2. \label{ineq:D21:2}
	\end{align}
	The first and second terms on the right hand side of \eqref{ineq:D21:2} are controllable since $j_0 \not= i$ and $i_1 \notin \{ i, l \}$. Meanwhile the higher order derivatives in the last term are integrable since $j_0 \not= i$ and $i_1 \not= i$. 
	
	{\myemph The case $(i,j)\in \mathcal{D}_{31}$. }
	
	Now we have $i_1 = i \not= j_0$ and $j_0 \not= j$. Let $l \in \{ 1, 2, 3 \} \setminus \{ j, j_0 \}$, then $\{ j, j_0, l \}$ is a permutation of $\{ 1, 2, 3 \}$. Use $f = \pt_i u_1$, $g = \pt_j b_1$, $h = \pt_i b_j$ in Lemma \ref{lem:case:1}, we have 
	\begin{align}
		& \Big| \int \pt_i b_j \pt_j b_1 \pt_i u_1 dx \Big| \nonumber \\
		& \le C\|\partial_i u_1\|_{L^2} \|\Lambda_j^\frac{1}{4}\partial_j b_1\|_{L^2}^\frac{1}{2} \|\Lambda_j^\frac{1}{4}\Lambda_{j_0}\partial_j b_1\|_{L^2}^\frac{1}{2} \|\Lambda_j^\frac{1}{4}\partial_i b_j\|_{L^2}^\frac{1}{2} \|(\Lambda_j^\frac{5}{4},\Lambda_l^\frac{5}{4})\partial_i b_j\|_{L^2}^\frac{1}{2} \nonumber \\
		& \le \frac{\eta}{100}\|\Lambda_j^\frac{5}{4}\partial_{j_0} b_1\|_{L^2}^2 + \frac{\eta}{100}\|(\Lambda_j^\frac{5}{4},\Lambda_l^\frac{5}{4})\partial_i b_j
		\|_{L^2}^2+ C\|\Lambda_j^\frac{5}{4} b_1\|_{L^2}\|(\Lambda_j^\frac{5}{4},\Lambda_i^\frac{5}{4}) b_j\|_{L^2}\|\partial_i u_1\|_{L^2}^2. \label{ineq:D31}
	\end{align}
	The first and second terms on the right hand side of \eqref{ineq:D31} are controllable since $j_0 \notin \{ j, l \}$. Meanwhile the higher order derivatives in the last term are integrable since $j_0 \notin \{ j, i \}$. 
	
	{\myemph The case $(i,j)\in \mathcal{D}_{42}$. }
	
	Now we have $i_1 = i \not= j_0 = j$. Since $i_1 \not= j_0$, let $l \in \{ 1, 2, 3 \}\setminus \{ i_1, j_0 \}$, then $\{ i_1, j_0, l \}$ is a permutation of $\{ 1, 2, 3 \}$. 
	
	Integration by parts yields
	\[
		\int \pt_i b_j \pt_j b_1 \pt_i u_1 dx = - \int u_1 \pt_i \pt_i b_j \pt_j b_1 dx - \int u_1 \pt_i \pt_j b_1 \pt_i b_j dx =: I_a + I_b. 
	\]
	
	Use $f = \pt_i \pt_i b_j$, $g = u_1$, $h = \pt_j b_1$ in Lemma \ref{lem:case:3}, we have
	\begin{align}
		|I_a| & \le C\|\Lambda_{i_1}^\frac{1}{4}\pt_i\pt_ib_j\|_{L^2} \|\Lambda_{j_0}^\frac{5}{4}u_1\|_{L^2}^\frac{2}{5} \| u_1\|_{L^2}^\frac{3}{5} \|\pt_jb_1\|_{L^2}^\frac{2}{5} \|(\Lambda_{i_1}^\frac{5}{4},\Lambda_l^\frac{5}{4}) \pt_jb_1\|_{L^2}^\frac{3}{5} \nonumber \\
		& \le \frac{\eta}{100}\|(\Lambda_{i_1}^\frac{5}{4},\Lambda_i^\frac{5}{4})\pt_ib_j\|_{L^2}^2 + \frac{\eta}{100}\|(\Lambda_{i_1}^\frac{5}{4},\Lambda_l^\frac{5}{4})\pt_jb_1\|_{L^2}^2 
		+ C\|\Lambda_{j_0}^\frac{5}{4} u_1\|_{L^2}^2 \|u_1\|_{L^2}^3 \|\pt_j b_1\|_{L^2}^2. \label{ineq:D42:1}
	\end{align}
	The first and second terms on the right hand side of \eqref{ineq:D42:1} are controllable since $j_0 \notin \{i,i_1\}$ and $j_0 \notin \{i_1,l\}$. Meanwhile the higher order derivative in the last term is integrable since $i_1 \not= j_0$.
	
	Use $f=\pt_i\pt_jb_1$, $g=\pt_ib_j$, $h=u_1$ in Lemma \ref{lem:case:3}, we have
	\begin{align}
		|I_b| & \le C\|\Lambda_l^\frac{1}{4}\pt_i\pt_jb_1\|_{L^2} \|\Lambda_{i_1}^\frac{5}{4}\pt_ib_j\|_{L^2}^\frac{2}{5} \| \pt_ib_j\|_{L^2}^\frac{3}{5} \|u_1\|_{L^2}^\frac{2}{5} \|(\Lambda_{j_0}^\frac{5}{4},\Lambda_l^\frac{5}{4}) u_1\|_{L^2}^\frac{3}{5} \nonumber \\
		& \le \frac{\eta}{100}\|(\Lambda_l^\frac{5}{4},\Lambda_i^\frac{5}{4})\pt_jb_1\|_{L^2}^2 + \frac{\eta}{100}\|\Lambda_{i_1}^\frac{5}{4}\pt_ib_j\|_{L^2}^2 
		+ C\|(\Lambda_{j_0}^\frac{5}{4},\Lambda_l^\frac{5}{4}) u_1\|_{L^2}^2 \|u_1\|_{L^2}^\frac{4}{3} \|\pt_ib_j\|_{L^2}^2. \label{ineq:D42:2}
	\end{align}
	The first and second terms on the right hand side of \eqref{ineq:D42:2} are controllable since $j_0 \notin\{l,i\}$ and $j_0 \not= i_1$. Meanwhile the higher order derivatives in the last term are integrable since $i_1\notin\{j_0, l\}$.
	
	{\myemph The case $(i,j)\in \mathcal{D}_{22}$. }
	
	Now we have $j_0 = i  = j = 1 \not= i_1$. We obtain since ${\bm b}$ is divergence free that
	\[
		\int \pt_1 b_1 \pt_1 b_1 \pt_1 u_1 dx = - \int \pt_1 b_1 ( \pt_2b_2 + \pt_3 b_3 ) \pt_1 u_1 dx =: I_a + I_b.
	\]
	
	Let $l\in\{1,2,3\}\setminus \{1, i_1\}$, then $\{1, i_1, l\}$ is a permutation of $\{1, 2, 3\}$. Use $f = \pt_1 b_1$, $g = \pt_2 b_2$, $h = \pt_1 u_1$ in Lemma \ref{lem:case:1}, we have 
	\begin{align}
		|I_a| & \le C\|\partial_1 b_1\|_{L^2} \|\Lambda_l^\frac{1}{4}\partial_2 b_2\|_{L^2}^\frac{1}{2} \|\Lambda_l^\frac{1}{4}\Lambda_{i_1}\partial_2 b_2\|_{L^2}^\frac{1}{2} \|\Lambda_l^\frac{1}{4}\partial_1 u_1\|_{L^2}^\frac{1}{2} \|(\Lambda_l^\frac{5}{4},\Lambda_{1}^\frac{5}{4})\partial_1 u_1\|_{L^2}^\frac{1}{2} \nonumber \\
		& \le \frac{\eta}{100}\|(\Lambda_l^\frac{5}{4},\Lambda_{i_1}^\frac{5}{4})\partial_2 b_2\|_{L^2}^2 + \frac{\nu}{100}\|(\Lambda_l^\frac{5}{4},\Lambda_{1}^\frac{5}{4})\partial_1 u_1 \|_{L^2}^2 \nonumber \\
		& + C \|(\Lambda_l^\frac{5}{4},\Lambda_2^\frac{5}{4}) b_2\|_{L^2}\|(\Lambda_l^\frac{5}{4},\Lambda_1^\frac{5}{4}) u_1\|_{L^2}\|\partial_1 b_1\|_{L^2}^2. \label{ineq:D22}
	\end{align}
	The first and second terms on the right hand side of \eqref{ineq:D22} are controllable since $j_0 = 1 \notin \{ l, i_1 \}$ and $i_1 \notin \{ 1, l \}$. Meanwhile the higher order derivatives in the last term are integrable since $j_0 = 1 \notin \{2,l\}$ and $i_1\notin\{1, l\}$. The estimates of $I_c$ are similar to those of $I_a$. We omit details for brevity. 
	
	{\myemph The case $(i,j)\in \mathcal{D}_{23}$. }
	
	Now we have $j_0 = i = j \not= i_1 = 1$. Let $l \in \{1,2,3\}\setminus \{1,i\}$, then $\{1,i,l\}$ is a permutation of $\{1,2,3\}$. Due to the solenoidal condition, we have 
	\[
		\int \pt_i b_i \pt_i b_1 \pt_i u_1 dx = - \int ( \pt_1 b_1 + \pt_l b_l ) \pt_i b_1 \pt_i u_1 dx =: I_a + I_b. 
	\]
	
	Use $f = \pt_i b_1$, $g = \pt_1 b_1$, $h = \pt_i u_1$ in Lemma \ref{lem:case:1}, we have
	\begin{align}
		| I_a | & \le C \| \pt_i b_1 \|_{L^2} \| \Lambda_l^\frac{1}{4}\pt_1 b_1 \|_{L^2}^\frac{1}{2} \| \Lambda_l^\frac{1}{4}\Lambda_1 \pt_1 b_1 \|_{L^2}^\frac{1}{2} \| \Lambda_l^\frac{1}{4} \pt_i u_1 \|_{L^2}^\frac{1}{2} \| ( \Lambda_l^\frac{5}{4},\Lambda_i^\frac{5}{4} ) \pt_i u_1 \|_{L^2}^\frac{1}{2} \nonumber \\
		& \le \frac{\eta}{100} \| ( \Lambda_l^\frac{5}{4}, \Lambda_1^\frac{5}{4} ) \pt_1 b_1 \|_{L^2}^2 + \frac{\nu}{100} \| ( \Lambda_l^\frac{5}{4}, \Lambda_i^\frac{5}{4} ) \pt_i u_1 \|_{L^2}^2 \nonumber \\
		& + C \| ( \Lambda_l^\frac{5}{4}, \Lambda_1^\frac{5}{4} ) b_1 \|_{L^2} \| ( \Lambda_l^\frac{5}{4}, \Lambda_i^\frac{5}{4} ) u_1 \|_{L^2} \| \pt_i b_1 \|_{L^2}^2. \label{ineq:D23:1}
	\end{align}
	The first and second terms on the right hand side of \eqref{ineq:D23:1} are controllable since $j_0 = i \notin \{ 1, l \}$ and $i_1 = 1 \notin \{ i, l \}$. And the higher order derivatives in the last term are integrable. 
	
	Since $j_0 = i \not= l$. Use $f = \pt_i b_1$, $g=\pt_lb_l$, $h=\pt_iu_1$ in Lemma \ref{lem:case:1}, we have
	\begin{align}
		|I_b| & \le C \|\pt_i b_1\|_{L^2} \|\Lambda_l^\frac{1}{4}\pt_l b_l\|_{L^2}^\frac{1}{2} \|\Lambda_l^\frac{1}{4}\Lambda_1\pt_l b_l\|_{L^2}^\frac{1}{2} \|\Lambda_l^\frac{1}{4}\pt_iu_1\|_{L^2}^\frac{1}{2} \|(\Lambda_l^\frac{5}{4},\Lambda_i^\frac{5}{4})\pt_i u_1\|_{L^2}^\frac{1}{2} \nonumber \\
		& \le \frac{\eta}{100}\|\Lambda_l^\frac{5}{4}\pt_1b_l\|_{L^2}^2 +\frac{\nu}{100}\|(\Lambda_l^\frac{5}{4},\Lambda_i^\frac{5}{4})\pt_iu_1\|_{L^2}^2 + C \| \Lambda_l^\frac{5}{4} b_l\|_{L^2}\|(\Lambda_l^\frac{5}{4},\Lambda_i^\frac{5}{4}) u_1\|_{L^2}\|\pt_ib_1\|_{L^2}^2. \label{ineq:D23:2}
	\end{align}
	The first and second terms on the right hand side of \eqref{ineq:D23:2} are controllable since $j_0 \not= l$ and $i_1 = 1 \notin \{ i, l \}$. Meanwhile the higher order derivatives in the last term are integrable since $j_0 \not= l$ and $i_1 \notin \{ i, l \}$. 
	
	{\myemph The case $(i, j)\in \mathcal{D}_{24}$. }
	
	Now we have $j_0 = i = j \not= i_1$, $i\not= 1$ and $i_1 \not= 1$. Then $\{ 1, i, i_1 \}$ is a permutation of $\{1,2,3\}$. Due to the solenoidal condition, we have 
	\[
		\int \pt_i b_i \pt_i b_1 \pt_i u_1 dx = - \int \pt_1 b_1 \pt_i b_1 \pt_i u_1 dx - \int \pt_{i_1} b_{i_1} \pt_i b_1 \pt_i u_1 dx. 
	\]
	The estimates of the above terms are similar to those in case $\mathcal{D}_{22}$. We omit the details for brevity. 
	
	{\myemph The case $(i, j) \in \mathcal{D}_{32}$. } 
	
	Now we have $j_0 = i_1 = i \not= j$. According to the definition of $\mathcal{J}$, it must hold that $j_0 = i_1 = 1$, then $j \not= 1$. Let $l \in \{ 1, 2, 3 \} \setminus \{ 1, j \}$, then $\{1,j,l\}$ is a permutation of $\{1,2,3\}$. Due to the solenoidal condition, we have 
	\[
		\int \pt_1 b_j \pt_j b_1 \pt_1 u_1 dx = - \int \pt_1 b_j \pt_j b_1 \pt_j u_j dx - \int \pt_1 b_j \pt_j b_1 \pt_l u_l dx =: I_a + I_b. 
	\]
	
	Integration by parts yields 
	\[
		I_a = \int u_j \pt_j \pt_1 b_j \pt_j b_1 dx + \int u_j \pt_j \pt_j b_1 \pt_1 b_j dx =: I_c + I_d. 
	\]
	
	Use $f = \pt_j \pt_1 b_j$, $g = \pt_j b_1$, $h = u_j$ in Lemma \ref{lem:case:2}, we have
	\begin{align}
		| I_c | & \le C \| \Lambda_l^\frac{1}{4} \pt_j \pt_1 b_j \|_{L^2} \| \Lambda_j^\frac{1}{4} \pt_j b_1 \|_{L^2}^\frac{1}{2} \| \Lambda_j^\frac{1}{4} \Lambda_1 \pt_j b_1 \|_{L^2}^\frac{1}{2} \| u_j \|_{L^2}^\frac{1}{2} \| ( \pt_l, \pt_j ) u_j \|_{L^2}^\frac{1}{2} \nonumber \\
		& \le \frac{\eta}{100} \| \Lambda_j^\frac{5}{4} \partial_1 b_1 \|_{L^2}^2 + \frac{\eta}{100} \| ( \Lambda_l^\frac{5}{4}, \Lambda_j^\frac{5}{4} ) \pt_1 b_j \|_{L^2}^2 + C \| u_j \|_{L^2}^2 \| \Lambda_j^\frac{5}{4} b_1 \|_{L^2}^2 \| ( \pt_l,\pt_j ) u_j \|_{L^2}^2. \label{ineq:D32:1}
	\end{align}
	The first and second terms on the right hand side of \eqref{ineq:D32:1} are controllable since $j_0 \not= j$ and $j_0 = 1 \notin\{l,j\}$. Meanwhile the higher order derivative in the last term is integrable since $j_0 \not= j$. 
	
	If $i_j=1$, then $i_j=j_0=1$ with $j \not= 1$,  $(i_1,i_2,i_3,j_0)\in \mathcal{J}$. Hence $i_j \not= 1$. Use $f=\pt_j\pt_jb_1$, $g=u_j$, $h=\pt_1b_j$ in Lemma \ref{lem:case:3}, we have
	\begin{align}
		|I_d| & \le C \| \Lambda_j^\frac{1}{4} \pt_j \pt_j b_1 \|_{L^2} \| \Lambda_1^\frac{5}{4} u_j \|_{L^2}^\frac{2}{5} \| u_j \|_{L^2}^\frac{3}{5} \| \pt_1 b_j \|_{L^2}^\frac{2}{5} \|( \Lambda_l^\frac{5}{4}, \Lambda_j^\frac{5}{4} ) \pt_1 b_j \|_{L^2}^\frac{3}{5} \nonumber \\
		& \le \frac{\eta}{100} \| \Lambda_j^\frac{5}{4} \pt_j b_1 \|_{L^2}^2 + \frac{\eta}{100} \| ( \Lambda_l^\frac{5}{4}, \Lambda_j^\frac{5}{4} ) \pt_1 b_j \|_{L^2}^2 
		+ C \| \Lambda_1^\frac{5}{4} u_j \|_{L^2}^2 \| u_j \|_{L^2}^3 \| \pt_1 b_j \|_{L^2}^2. \label{ineq:D32:2}
	\end{align}
	The first and second terms on the right hand side of \eqref{ineq:D32:2} are controllable since $j_0 = 1 \notin \{ l, j \}$. Meanwhile the higher order derivative in the last term is integrable since $i_j \not= 1$. 
	
	Integration by parts yields
	\[
		I_b = \int u_l \pt_l \pt_1 b_j \pt_j b_1 dx + \int u_l \pt_l \pt_j b_1 \pt_1 b_j dx =: I_e + I_f. 
	\]
	
	Use $f = \pt_l \pt_1 b_j$, $g = \pt_j b_1$, $h = u_l$ in Lemma \ref{lem:case:2}, we have 
	\begin{align}
		| I_e | & \le C \| \Lambda_l^\frac{1}{4} \pt_l \pt_1 b_j \|_{L^2} \| \Lambda_j^\frac{1}{4} \pt_j b_1 \|_{L^2}^\frac{1}{2} \| \Lambda_j^\frac{1}{4} \Lambda_1 \pt_j b_1 \|_{L^2}^\frac{1}{2} \| u_l \|_{L^2}^\frac{1}{2} \| ( \pt_l, \pt_j ) u_l \|_{L^2}^\frac{1}{2} \nonumber \\
		& \le \frac{\eta}{100} \| \Lambda_j^\frac{5}{4} \partial_1 b_1 \|_{L^2}^2 + \frac{\eta}{100} \| \Lambda_l^\frac{5}{4} \pt_1 b_j \|_{L^2}^2 + C \| u_l \|_{L^2}^2 \| \Lambda_j^\frac{5}{4} b_1 \|_{L^2}^2 \| ( \pt_l, \pt_j ) u_l \|_{L^2}^2. \label{ineq:D32:3}
	\end{align}
	The first and second terms on the right hand side of \eqref{ineq:D32:3} are controllable since $j_0 \notin \{j,l\}$. Meanwhile the higher order derivative in the last term is integrable since $j_0 \not= j$. 
	
	If $i_l=1$, then $i_l=j_0=1$ with $l \not= 1$, $(i_1,i_2,i_3,j_0)\in \mathcal{J}$. Hence $i_l \not= 1$. Use $f=\pt_l\pt_jb_1$, $g=u_l$, $h=\pt_1b_j$ in Lemma \ref{lem:case:3}, we have
	\begin{align}
		|I_f| & \le C\|\Lambda_j^\frac{1}{4}\pt_l\pt_jb_1\|_{L^2} \|\Lambda_1^\frac{5}{4}u_l\|_{L^2}^\frac{2}{5} \| u_l\|_{L^2}^\frac{3}{5} \|\pt_1b_j\|_{L^2}^\frac{2}{5} \|(\Lambda_l^\frac{5}{4},\Lambda_j^\frac{5}{4}) \pt_1b_j\|_{L^2}^\frac{3}{5} \nonumber \\
		& \le \frac{\eta}{100}\|(\Lambda_l^\frac{5}{4},\Lambda_j^\frac{5}{4})\pt_j b_1\|_{L^2}^2 + \frac{\eta}{100}\|(\Lambda_l^\frac{5}{4},\Lambda_j^\frac{5}{4})\pt_1b_j\|_{L^2}^2 
		+ C\|\Lambda_1^\frac{5}{4} u_l\|_{L^2}^2 \| u_l\|_{L^2}^3 \|\pt_1 b_j\|_{L^2}^2. \label{ineq:D32:4}
	\end{align}
	The first and second terms on the right hand side of \eqref{ineq:D32:4} are controllable since $j_0 = 1 \notin\{l,j\}$. The higher order derivative in the last term is integrable since $i_l \not= 1$. 
	
	{\myemph The case $(i,j)\in \mathcal{D}_{41}$. }
	
	Now we have $i_1 = j_0 = i = j$. According to the set $\mathcal{J}$, we have $i_1 = j_0 = i = j = 1$. Due to the solenoidal condition, we have 
	\[
		\int ( \pt_1 b_1 )^2 \pt_1 u_1 dx = - \int(\pt_1b_1)^2(\pt_2u_2+\pt_3u_3)dx
		=: I_a + I_b.
	\]
	
	If $i_2 = 1$, then $i_2 = j_0 = 1$, $( i_1, i_2, i_3, j_0 ) \in \mathcal{J}$. Hence $i_2 \not= 1$. Integration by parts yields 
	\[
		I_a=2\int u_2\pt_2\pt_1b_1\pt_1b_1dx.
	\]
	Use $f=\pt_2\pt_1b_1$, $g=u_2$, $h=\pt_1b_1$ in Lemma \ref{lem:case:3}, we have
	\begin{align}
		|I_a| & \le C\|\Lambda_2^\frac{1}{4}\pt_2\pt_1b_1\|_{L^2} \|\Lambda_1^\frac{5}{4}u_2\|_{L^2}^\frac{2}{5} \| u_2\|_{L^2}^\frac{3}{5} \|\pt_1b_1\|_{L^2}^\frac{2}{5} \|(\Lambda_2^\frac{5}{4},\Lambda_3^\frac{5}{4}) \pt_1b_1\|_{L^2}^\frac{3}{5} \nonumber \\
		& \le \frac{\eta}{100}\|\Lambda_2^\frac{5}{4}\pt_1b_1\|_{L^2}^2 + \frac{\eta}{100}\|(\Lambda_2^\frac{5}{4},\Lambda_3^\frac{5}{4})\pt_1b_1\|_{L^2}^2 
		+ C\|\Lambda_1^\frac{5}{4} u_2\|_{L^2}^2 \|u_2\|_{L^2}^3 \|\pt_1 b_1\|_{L^2}^2. \label{ineq:D41}
	\end{align}
	The first and second terms on the right hand side of \eqref{ineq:D41} are controllable since $j_0 = 1 \notin\{2,3\}$. Meanwhile the higher order derivative in the last term is integrable since $i_2\not=1$. The estimate of $I_b$ is similar to that of $I_a$. We omit the details for brevity. 
	
	\subsection{The estimates of $K_4$}
	We now estimate the terms in $K_4$,
	\begin{align*}
		K_4 = & - \int \partial _i u_j \cdot \partial _j b_k \cdot \partial _i b_k dx\\
		= & \sum\limits_{i,j = 1}^{\rm 3} \Big[ -\int \partial _i u_j\partial _j b_1\partial _i b_1 
		- \partial _i u_j\partial _j b_2\partial _i b_2 
		- \partial _i u_j\partial _j b_3\partial _i b_3 \Big] dx =: K_{41} + K_{42} + K_{43}. 
	\end{align*}
	We only estimate $K_{41}$ in the following context since $K_{42}$, $K_{43}$ can be argued similarly. Now we separate cases according to the indices $i$ and $j$. Define
	\[
	\begin{array}{l}
		\mathcal{E}_1 := \{(i,j)\in\{1,2,3\} \mid i \not= j_0, j \not= j_0 \}, \\
		\mathcal{E}_2 := \{(i,j)\in\{1,2,3\} \mid i = j_0, j \not= j_0 \}, \\
		\mathcal{E}_3 := \{(i,j)\in\{1,2,3\} \mid i \not= j_0, j = j_0 \} = \mathcal{E}_{31} \cup \mathcal{E}_{32}, \\
		\mathcal{E}_4 := \{(i,j)\in\{1,2,3\} \mid i = j_0, j = j_0 \} = \mathcal{E}_{41} \cup \mathcal{E}_{42}, 
	\end{array}
	\]
	where
	\[
	\begin{array}{ll}
		\mathcal{E}_{31} := \{(i,j) \in \mathcal{E}_3 \mid i_j \not= i \}, & \qquad \mathcal{E}_{32} := \{(i,j) \in \mathcal{E}_3 \mid i_j = i\}, \\
		\mathcal{E}_{41} := \{(i,j) \in \mathcal{E}_4 \mid j_0 \not= 1 \}, & \qquad \mathcal{E}_{42} := \{(i,j) \in \mathcal{E}_4 \mid j_0 = 1 \}. 
	\end{array}
	\]
	\setcounter{mystep}{0}
	{\myemph The case $(i,j)\in \mathcal{E}_1$.}
	
	Since $i \not=j_0$ and $j \not= j_0$, let $k\in \{1,2,3\}\setminus \{i, j_0\}$, then $\{i, j_0, k\}$ is a permutation of $\{1,2,3\}$. Use $f=\pt_iu_j$, $g=\pt_ib_1$, $h=\pt_jb_1$ in Lemma \ref{lem:case:1}, we have
	\begin{align}
		& |\int\pt_iu_j\pt_jb_1\pt_ib_1dx| \nonumber \\
		& \le C\|\pt_iu_j\|_{L^2} \|\Lambda_i^\frac{1}{4}\pt_i b_1\|_{L^2}^\frac{1}{2} \|\Lambda_i^\frac{1}{4} \Lambda_{j_0}\pt_i b_1\|_{L^2}^\frac{1}{2} \|\Lambda_i^\frac{1}{4}\pt_j b_1\|_{L^2}^\frac{1}{2} \|(\Lambda_i^\frac{5}{4},\Lambda_k^\frac{5}{4})\pt_j b_1\|_{L^2}^\frac{1}{2} \nonumber \\
		& \le \frac{\eta}{100}\|(\Lambda_i^\frac{5}{4},\Lambda_k^\frac{5}{4})\pt_j b_1\|_{L^2}^2 + \frac{\eta}{100}\| \Lambda_i^\frac{5}{4} \pt_{j_0} b_1\|_{L^2}^2+C\|(\Lambda_i^\frac{5}{4},\Lambda_j^\frac{5}{4}) b_1\|^2_{L^2}\|\pt_iu_j\|_{L^2}^2. \label{ineq:E1}
	\end{align}
	The first and second terms on the right hand side of \eqref{ineq:E1} are controllable since $j_0 \notin \{ i, k \}$. Meanwhile the higher order derivatives in the last term are integrable since $j_0 \notin \{ i, j \}$. 
	
	{\myemph The case $(i,j)\in \mathcal{E}_{2}$.}
	
	If $i_j=i$, then $i_j = j_0 \not= j$, and $(i_1,i_2,i_3,j_0)\in \mathcal{J}$. Hence $i_j \not= i$. Now we have $ i = j_0 \not= i_j$ and $j \not= j_0$, let $l\in \{1,2,3\}\setminus \{i,i_j\}$, then $\{i,i_j,l\}$ is a permutation of $\{1,2,3\}$. Use $f=\pt_ib_1$, $g=\pt_iu_j$, $h=\pt_jb_1$ in Lemma \ref{lem:case:1}, we have
	\begin{align}
		& |\int\pt_iu_j\pt_jb_1\pt_ib_1dx| \nonumber \\
		& \le C\|\pt_ib_1\|_{L^2} \|\Lambda_l^\frac{1}{4}\pt_i u_j\|_{L^2}^\frac{1}{2} \|\Lambda_l^\frac{1}{4} \Lambda_i\pt_i u_j\|_{L^2}^\frac{1}{2} \|\Lambda_l^\frac{1}{4}\pt_j b_1\|_{L^2}^\frac{1}{2} \|(\Lambda_l^\frac{5}{4},\Lambda_{i_j}^\frac{5}{4})\pt_j b_1\|_{L^2}^\frac{1}{2} \nonumber \\
		& \le \frac{\eta}{100}\|(\Lambda_l^\frac{5}{4},\Lambda_{i_j}^\frac{5}{4})\pt_j b_1\|_{L^2}^2 + \frac{\nu}{100}\| (\Lambda_l^\frac{5}{4},\Lambda_i^\frac{5}{4})\pt_i u_j\|_{L^2}^2 \nonumber \\
		& + C \|(\Lambda_l^\frac{5}{4},\Lambda_i^\frac{5}{4}) u_j\|_{L^2}\|(\Lambda_l^\frac{5}{4},\Lambda_j^\frac{5}{4}) b_1\|_{L^2}\|\pt_ib_1\|_{L^2}^2. \label{ineq:E2}
	\end{align}
	The first and second terms on the right hand side of \eqref{ineq:E2} are controllable since $j_0 = i \notin\{l,i_j\}$ and $i_j\notin\{l,i\}$. Meanwhile the higher order derivatives in the last term are integrable since $i_j\notin\{l,i\}$, $j_0 \notin \{ l, j \}$. 

	{\myemph The case $(i,j)\in \mathcal{E}_{31}$.}
	
	Since $i \not= j = j_0$ and $i \not= i_j$, let $l\in \{1,2,3\}\setminus \{i,j\}$, then $\{i,j,l\}$ is a permutation of $\{1,2,3\}$. Use $f=\pt_jb_1$, $g=\pt_iu_j$, $h=\pt_ib_1$ in Lemma \ref{lem:case:1}, we have
	\begin{align}
		& |\int\pt_iu_j\pt_jb_1\pt_ib_1dx| \nonumber \\
		& \le C\|\pt_jb_1\|_{L^2} \|\Lambda_i^\frac{1}{4}\pt_i u_j\|_{L^2}^\frac{1}{2} \|\Lambda_i^\frac{1}{4} \Lambda_j\pt_i u_j\|_{L^2}^\frac{1}{2} \|\Lambda_i^\frac{1}{4}\pt_i b_1\|_{L^2}^\frac{1}{2} \|(\Lambda_i^\frac{5}{4},\Lambda_l^\frac{5}{4})\pt_i b_1\|_{L^2}^\frac{1}{2} \nonumber \\
		& \le \frac{\eta}{100}\|(\Lambda_i^\frac{5}{4},\Lambda_l^\frac{5}{4})\pt_i b_1\|_{L^2}^2 + \frac{\nu}{100}\| \Lambda_i^\frac{5}{4} \pt_j u_j\|_{L^2}^2+C\|\Lambda_i^\frac{5}{4} b_1\|_{L^2}\|\Lambda_i^\frac{5}{4} u_j\|_{L^2}\|\pt_jb_1\|_{L^2}^2. \label{ineq:E31}
	\end{align}
	The first and second terms on the right hand side of \eqref{ineq:E31} are controllable since $j_0 = j \notin \{i,l\}$ and $i_j\not=i$. Meanwhile the higher order derivatives in the last term are integrable. 
	
	{\myemph The case $(i,j)\in \mathcal{E}_{32}$.}
	
	Since $i = i_j \not= j = j_0$, let $l\in \{1,2,3\} \setminus \{i,j\}$, then $\{i,j,l\}$ is a permutation of $\{1,2,3\}$. Integration by parts yields 
	\[
		-\int\pt_iu_j\pt_jb_1\pt_ib_1dx=\int u_j\pt_i\pt_jb_1\pt_ib_1dx+\int u_j\pt_i\pt_ib_1\pt_jb_1 dx=:I_a+I_b.
	\]
	
	Use $f=\pt_i\pt_jb_1$, $g=u_j$, $h=\pt_ib_1$ in Lemma \ref{lem:case:3}, we have
	\begin{align}
		|I_a| & \le C\|\Lambda_i^\frac{1}{4}\pt_i\pt_jb_1\|_{L^2} \|\Lambda_j^\frac{5}{4}u_j\|_{L^2}^\frac{2}{5} \| u_j\|_{L^2}^\frac{3}{5} \|\pt_ib_1\|_{L^2}^\frac{2}{5} \|(\Lambda_i^\frac{5}{4},\Lambda_l^\frac{5}{4}) \pt_ib_1\|_{L^2}^\frac{3}{5} \nonumber \\
		& \le \frac{\eta}{100}\|\Lambda_i^\frac{5}{4}\pt_jb_1\|_{L^2}^2 + \frac{\eta}{100}\|(\Lambda_i^\frac{5}{4},\Lambda_l^\frac{5}{4})\pt_ib_1\|_{L^2}^2 
		+ C\|\Lambda_j^\frac{5}{4} u_j\|_{L^2}^2 \|u_j\|_{L^2}^3 \|\pt_i b_1\|_{L^2}^2. \label{ineq:E32:1}
	\end{align}
	The first and second terms on the right hand side of \eqref{ineq:E32:1} are controllable since $j_0 = j \notin \{ i, l \}$. Meanwhile the higher order derivative in the last term is integrable since $i_j = i \not= j$. 
	
	Use $f=\pt_i\pt_ib_1$, $g=u_j$, $h=\pt_jb_1$ in Lemma \ref{lem:case:3}, we have
	\begin{align}
		|I_b| & \le C\|\Lambda_i^\frac{1}{4}\pt_i\pt_ib_1\|_{L^2} \|\Lambda_j^\frac{5}{4}u_j\|_{L^2}^\frac{2}{5} \| u_j\|_{L^2}^\frac{3}{5} \|\pt_jb_1\|_{L^2}^\frac{2}{5} \|(\Lambda_i^\frac{5}{4},\Lambda_l^\frac{5}{4}) \pt_jb_1\|_{L^2}^\frac{3}{5} \nonumber \\
		& \le \frac{\eta}{100}\|\Lambda_i^\frac{5}{4}\pt_ib_1\|_{L^2}^2 + \frac{\eta}{100}\|(\Lambda_i^\frac{5}{4},\Lambda_l^\frac{5}{4})\pt_jb_1\|_{L^2}^2 
		+ C\|\Lambda_j^\frac{5}{4} u_j\|_{L^2}^2 \|u_j\|_{L^2}^3 \|\pt_j b_1\|_{L^2}^2. \label{ineq:E32:2}
	\end{align}
	The first and second terms on the right hand side of \eqref{ineq:E32:2} are controllable since $j_0 = j \notin \{ i, l \}$. The higher order derivative in the last term is integrable since $i_j=i\not=j$. 
	
	{\myemph The case $(i,j)\in \mathcal{E}_{41}$.}
	
	Since $i = j = j_0 \not= 1$, let $l\in \{1,2,3\}\setminus \{1,i\}$, then $\{1,i,l\}$ is a permutation of $\{1,2,3\}$. We obtain since ${\bm u}$ is divergence free that
	\[
		-\int\pt_iu_i\pt_ib_1\pt_ib_1dx=\int\pt_1u_1\pt_ib_1\pt_ib_1dx+\int\pt_lu_l\pt_ib_1\pt_ib_1dx=:I_a+I_b.
	\]
	
	Integration by parts allows us to have that
	\[
		I_a = - 2 \int u_1\pt_1\pt_ib_1\pt_ib_1dx.
	\]
	If $i_1 = i$, then $i_1 = j_0 \not= 1$, $(i_1,i_2,i_3,j_0) \in \mathcal{J}$. Hence $i_1 \not= i$. Use $f = \pt_1 \pt_i b_1$, $g = u_1$, $h = \pt_i b_1$ in Lemma \ref{lem:case:3}, we have
	\begin{align}
		|I_a| & \le C\|\Lambda_1^\frac{1}{4}\pt_1\pt_ib_1\|_{L^2} \|\Lambda_i^\frac{5}{4}u_1\|_{L^2}^\frac{2}{5} \| u_1\|_{L^2}^\frac{3}{5} \|\pt_ib_1\|_{L^2}^\frac{2}{5} \|(\Lambda_1^\frac{5}{4},\Lambda_l^\frac{5}{4}) \pt_ib_1\|_{L^2}^\frac{3}{5} \nonumber \\
		& \le \frac{\eta}{100}\|\Lambda_1^\frac{5}{4}\pt_ib_1\|_{L^2}^2 + \frac{\eta}{100}\|(\Lambda_1^\frac{5}{4},\Lambda_l^\frac{5}{4})\pt_ib_1\|_{L^2}^2 
		+ C\|\Lambda_i^\frac{5}{4} u_1\|_{L^2}^2 \|u_1\|_{L^2}^3 \|\pt_i b_1\|_{L^2}^2. \label{ineq:E41:1}
	\end{align}
	The first and second terms on the right hand side of \eqref{ineq:E41:1} are controllable since $j_0 = i \notin \{1,l\}$. The higher order derivative in the last term is integrable since $i_1 \not= i$. 
	
	From integration by parts, we have that
	\[
		I_b = - 2 \int u_l \pt_l \pt_i b_1 \pt_i b_1 dx. 
	\]
	If $i_l = i$, then $i_l = j_0 \not= 1$, $(i_1,i_2,i_3,j_0) \in \mathcal{J}$. Hence $i_l \not= i$. Use $f=\pt_l\pt_ib_1$, $g=u_l$, $h=\pt_ib_1$ in Lemma \ref{lem:case:3}, we have
	\begin{align}
		|I_b| & \le C\|\Lambda_l^\frac{1}{4}\pt_l\pt_ib_1\|_{L^2} \|\Lambda_i^\frac{5}{4}u_l\|_{L^2}^\frac{2}{5} \| u_l\|_{L^2}^\frac{3}{5} \|\pt_ib_1\|_{L^2}^\frac{2}{5} \|(\Lambda_1^\frac{5}{4},\Lambda_l^\frac{5}{4}) \pt_ib_1\|_{L^2}^\frac{3}{5} \nonumber \\
		& \le \frac{\eta}{100}\|\Lambda_l^\frac{5}{4}\pt_ib_1\|_{L^2}^2 + \frac{\eta}{100}\|(\Lambda_1^\frac{5}{4},\Lambda_l^\frac{5}{4})\pt_ib_1\|_{L^2}^2 
		+ C\|\Lambda_i^\frac{5}{4} u_l\|_{L^2}^2 \|u_l\|_{L^2}^3 \|\pt_i b_1\|_{L^2}^2. \label{ineq:E41:2}
	\end{align}
	The first and second terms on the right hand side of \eqref{ineq:E41:2} are controllable since $j_0 = i \notin \{ 1, l \}$. The higher order derivative in the last term is integrable since $i_l \not= i$. 
	
	{\myemph The case $(i,j)\in \mathcal{E}_{42}$.}
	
	Now we have $i = j = j_0 = 1$. We obtain from the solenoidal condition that
	\[
		-\int\pt_1u_1\pt_1b_1\pt_1b_1dx = \int\pt_2u_2 (\pt_1b_1)^2dx + \int \pt_3u_3 (\pt_1b_1)^2 dx =: I_a+I_b.
	\]
	The estimates of $I_b$ are similar to those of $I_a$. Here we only derive estimates for $I_a$. 
	
	If $i_2 = 1$, then $i_2 = j_0 = 1$, $(i_1,i_2,i_3,j_0) \in \mathcal{J}$. Hence $i_2 \not= 1$. Integration by parts yields that
	\[
		I_a=-2\int u_2\pt_1\pt_2b_1\pt_1b_1dx.
	\]
	Use $f=\pt_1\pt_2b_1$, $g=u_2$, $h=\pt_1b_1$ in Lemma \ref{lem:case:3}, we have
	\begin{align}
		|I_a| & \le C\|\Lambda_2^\frac{1}{4}\pt_1\pt_2b_1\|_{L^2} \|\Lambda_1^\frac{5}{4}u_2\|_{L^2}^\frac{2}{5} \| u_2\|_{L^2}^\frac{3}{5} \|\pt_1b_1\|_{L^2}^\frac{2}{5} \|(\Lambda_2^\frac{5}{4},\Lambda_3^\frac{5}{4}) \pt_1b_1\|_{L^2}^\frac{3}{5} \nonumber \\
		& \le \frac{\eta}{100}\|\Lambda_2^\frac{5}{4}\pt_1b_1\|_{L^2}^2 + \frac{\eta}{100}\|(\Lambda_2^\frac{5}{4},\Lambda_3^\frac{5}{4})\pt_1b_1\|_{L^2}^2 
		+ C\|\Lambda_1^\frac{5}{4} u_2\|_{L^2}^2 \|u_2\|_{L^2}^3 \|\pt_1 b_1\|_{L^2}^2. \label{ineq:E42}
	\end{align}
	The first and second terms on the right hand side of \eqref{ineq:E42} are controllable since $j_0 = 1 \notin \{2,3\}$. The higher order derivative in the last term is integrable since $i_2 \not= 1$. We first note that $i_3\not=1$, then the estimate of $I_b$ can be obtained similarly to those of $I_a$. 
	
	Now we have finished all required estimates for the proof of Theorem \ref{thm:existence:MHD}. 
	
	\vspace{10pt}
	
	\noindent
{\bf Continuum of the proof of Theorem \ref{thm:existence:MHD}.} 
	Combining all the above estimates together, we have that
	\begin{align*}
		& \frac{d}{dt}\| (\nabla {\bm u}, \nabla {\bm b}) \|_{L^2}^2+ \nu \sum_{j\not=i_1}\| \Lambda_j ^\frac{5}{4} \nabla u_1 \|_{L^2}^2 + \nu \sum_{j\not=i_2}\| \Lambda_j ^\frac{5}{4} \nabla u_2 \|_{L^2}^2 \\
		& + \nu \sum_{j\not=i_3}\| \Lambda_j ^\frac{5}{4} \nabla u_3 \|_{L^2}^2 + \eta \sum_{k=1}^{3} \sum_{j \not= j_0}\| \Lambda_j ^\frac{5}{4} \nabla b_k \|_{L^2}^2 \\
		& \le C\big( \sum_{j\not=i_1}\|\Lambda_j ^\frac{5}{4} u_1+\sum_{j\not=i_2}\|\Lambda_j ^\frac{5}{4} u_2 \|_{L^2}^2+\sum_{j\not=i_3}\|\Lambda_j ^\frac{5}{4} u_3\|_{L^2}^2 + \sum_{k=1}^{3} \sum_{j \not= j_0} \| \Lambda_j^\frac{5}{4} b_k \|_{L^2}^2 \big) \| ( \nabla {\bm u}, \nabla {\bm b} ) \|_{L^2}^2, 
	\end{align*}
	where $C$ may depend on $\| ({\bm u}_0, {\bm b}_0) \|_{L^2}^2$. The Gr\"{o}nwall's inequality then implies the desired uniform global $H^1$-bound as stated in Theorem \ref{thm:existence:MHD}.
\hfill$\square$

\section{Uniqueness for the MHD equations}

{\bf Proof of Theorem \ref{thm:uniqueness:MHD}.}
	It is sufficient to estimate the $L^2$-difference of two solutions
	\[
		(\widetilde{\bm{u}}, \widetilde{\bm{b}}) = (\bm{u}^{(1)}, \bm{b}^{(1)}) - (\bm{u}^{(2)}, \bm{b}^{(2)})
	\]
	which satisfies
	\begin{equation}\label{eq:MHD:different}
		\left\{ 
		\begin{aligned}
			& \partial _t \tilde{u}_1 + \boldsymbol{u}^{(1)} \cdot \nabla \tilde{u}_1+\widetilde{{\bm u}}\cdot \nabla u_1^{(2)} +\nu (\Lambda_1 ^\frac{5}{2} \tilde{u}_1+\Lambda _2^\frac{5}{2} \tilde{u}_1+\Lambda _3^\frac{5}{2} \tilde{u}_1-\Lambda _{i_1}^\frac{5}{2} \tilde{u}_1) + \partial _1 \tilde{p} 
			= {\bm b}^{(1)}\cdot\nabla \tilde{b}_1+\widetilde{{\bm b}}\cdot \nabla b_1^{(2)}, \\
			& \partial _t \tilde{u}_2 + \boldsymbol{u}^{(1)} \cdot \nabla \tilde{u}_2+\widetilde{{\bm u}}\cdot \nabla u_2^{(2)} +\nu (\Lambda_1 ^\frac{5}{2} \tilde{u}_2+\Lambda _2^\frac{5}{2} \tilde{u}_2+\Lambda _3^\frac{5}{2} \tilde{u}_2-\Lambda _{i_2}^\frac{5}{2} \tilde{u}_2) + \partial _2 \tilde{p} 
			= {\bm b}^{(1)}\cdot\nabla \tilde{b}_2+\widetilde{{\bm b}}\cdot \nabla b_2^{(2)}, \\
			& \partial _t \tilde{u}_3 + \boldsymbol{u}^{(1)} \cdot \nabla \tilde{u}_3+\widetilde{{\bm u}}\cdot \nabla u_3^{(2)}+\nu (\Lambda_1 ^\frac{5}{2} \tilde{u}_3+\Lambda _2^\frac{5}{2} \tilde{u}_3+\Lambda _3^\frac{5}{2} \tilde{u}_3-\Lambda _{i_3}^\frac{5}{2} \tilde{u}_3) + \partial _3 \tilde{p} 
			= {\bm b}^{(1)}\cdot\nabla \tilde{b}_3+\widetilde{{\bm b}}\cdot \nabla b_3^{(2)}, \\
			& \partial_t \tilde{\bm b} + {\bm u}^{(1)} \cdot \nabla \tilde{\bm b} +\widetilde{{\bm u}}\cdot \nabla {\bm b}^{(2)} + \eta (\Lambda_1 ^\frac{5}{2} \tilde{\bm b} + \Lambda _2^\frac{5}{2} \tilde{\bm b} + \Lambda _3^\frac{5}{2} \tilde{\bm b} - \Lambda _{j_0}^\frac{5}{2} \tilde{\bm b}) 
			= {\bm b}^{(1)} \cdot \nabla \tilde{\bm u} + \widetilde{{\bm b}}\cdot \nabla {\bm u}^{(2)}, \\
			& \dive \widetilde{\bm u} = 0, \qquad \dive \widetilde{\bm b} = 0,\\
			& \widetilde{\bm u}( x,0 ) = 0, \qquad \widetilde{\bm b}( x,0 ) = 0.
		\end{aligned} 
		\right.
	\end{equation}
	Here $ \tilde{p} := p^{( 1)} - p^{( 2 )} $ denotes the difference between the corresponding pressures. Taking inner product of \eqref{eq:MHD:different} and $(\widetilde{\bm{u}}, \widetilde{\bm{b}})$, then integrating by parts, we have
	\begin{align*}
		& \frac{1}{2} \frac{d}{dt}\| (\widetilde{\bm{u}}, \widetilde{\bm{b}}) \|_{L^2}^2 + \nu \sum_{j\not=i_1}\| \Lambda_j ^\frac{5}{4} \tilde{u}_1 \|_{L^2}^2 + \nu \sum_{j\not=i_2}\| \Lambda_j ^\frac{5}{4} \tilde{u}_2 \|_{L^2}^2 \\
		& + \nu \sum_{j\not=i_3}\| \Lambda_j ^\frac{5}{4} \tilde{u}_3 \|_{L^2}^2 
		+ \eta \sum_{k=1}^{3} \sum_{j \not= j_0} \| \Lambda_j ^\frac{5}{4} \tilde{b}_k \|_{L^2}^2 \\
		& = - \int \Big[( \widetilde{{\bm u}}\cdot \nabla {\bm u}^{(2)} ) \cdot \widetilde{{\bm u}} - ( \widetilde{{\bm b}}\cdot \nabla {\bm b}^{(2)} ) \cdot \widetilde{{\bm u}} +  ( \widetilde{{\bm u}}\cdot \nabla {\bm b}^{(2)} ) \cdot \widetilde{{\bm b}} - ( \widetilde{{\bm b}}\cdot \nabla {\bm u}^{(2)} ) \cdot \widetilde{{\bm b}}\Big] dx\\
		& =: L_1 + L_2 + L_3 + L_4. 
	\end{align*}
	The estimate of $L_1$ has been calculated in the section 4, which is controllable after imposing conditions on ${\bm u}^{(2)}$. The detailed calculation is omitted here. 
	
	\subsection{The estimates of $L_2$}
	
	The items in $L_2$ can be expanded as
	\begin{align*}
		L_2 & = \int \tilde{b}_k \partial_k b_m^{(2)} \tilde{u}_m dx 
		=\sum\limits_{k = 1}^{\rm 3} \Big[\int \tilde{b}_k \partial_k b_1^{(2)} \tilde{u}_1  + \tilde{b}_k \partial_k b_2^{(2)} \tilde{u}_2  + \tilde{b}_k \partial_k b_3^{(2)} \tilde{u}_3\Big] dx \\
		& =: L_{21} + L_{22} + L_{23}.
	\end{align*}
	We only estimate $L_{21}$ in the following context since $L_{22}$, $L_{23}$ can be argued similarly. Now we separate cases according to the indices $k$. Define
	\[
	\begin{array}{l}
		\mathcal{F}_1 := \{k\in\{1,2,3\} \mid i_1 \not= j_0 \}, \\
		\mathcal{F}_2 := \{k\in\{1,2,3\} \mid i_1 = j_0, k \not= 1 \}, \\
		\mathcal{F}_3 := \{k\in\{1,2,3\} \mid i_1 = j_0, k = 1 \}. 
	\end{array}
	\]
	
	\setcounter{mystep}{0}	
	{\myemph The case $k\in \mathcal{F}_{1}$. }
	
	Now we have $j_0 \not= i_1$. Let $l \in \{ 1, 2, 3 \} \setminus \{ i_1, j_0 \}$, then $\{ i_1, j_0, l \}$ is a permutation of $\{1,2,3\}$. Use $f=\tilde{b}_k$, $g=\partial_k b_1^{(2)}$, $h=\tilde{u}_1$ in Lemma \ref{lem:case:4}, we have
	\begin{align}
		& |\int \tilde{b}_k \partial_k b_1^{(2)} \tilde{u}_1 dx | \nonumber \\
		& \le C\|\tilde{b}_k\|_{L^2}^\frac{3}{5} \|\Lambda_l^\frac{5}{4}\tilde{b}_k\|_{L^2}^\frac{2}{5} \|\partial_k b_1^{(2)}\|_{L^2}^\frac{3}{5} \|\Lambda_{i_1}^\frac{5}{4}\partial_k b_1^{(2)}\|_{L^2}^\frac{2}{5} \|\tilde{u}_1\|_{L^2}^\frac{3}{5} \|\Lambda_{ j_0 }^\frac{5}{4}\tilde{u}_1\|_{L^2}^\frac{2}{5} \nonumber \\
		& \le \frac{\nu}{100}\|\Lambda_{ j_0 }^\frac{5}{4}\tilde{u}_1\|_{L^2}^2 + \frac{\eta}{100}\|\Lambda_l^\frac{5}{4}\tilde{b}_k\|_{L^2}^2 + C \| \tilde{b}_k \|_{L^2} \| \tilde{u}_1 \|_{L^2} \| \partial_k b_1^{(2)} \|_{L^2} \|\Lambda_{i_1}^\frac{5}{4}\partial_k b_1^{(2)}\|_{L^2}^\frac{2}{3}. \label{ineq:F1}
	\end{align}
	The first and second terms on the right hand side of \eqref{ineq:F1} are controllable, and the higher order derivative in the last term is integrable on any finite time interval since $i_1 \not= j_0$ and $j_0 \not= l$. 
	
	{\myemph The case $k\in \mathcal{F}_{2}$. }
	
	Now we have $i_1 = j_0, k \not= 1$. According to the definition of $\mathcal{J}$, we must have $i_1 = j_0 = 1$. Let $l\in\{1,2,3\}\setminus \{1,k\}$, then $\{1,k,l\}$ is a permutation of $\{1,2,3\}$. Use $f=\tilde{b}_k$, $g=\pt_kb_1^{(2)}$, $h=\tilde{u}_1$ in Lemma \ref{lem:case:1}, we have 
	\[
		| \int \tilde{b}_k \pt_k b_1^{(2)} \tilde{u}_1 dx | \le C \| \tilde{b}_k \|_{L^2} \| \Lambda_k^\frac{1}{4} \pt_k b_1^{(2)} \|_{L^2}^\frac{1}{2} \| \Lambda_k^\frac{1}{4} \Lambda_1 \pt_k b_1^{(2)} \|_{L^2}^\frac{1}{2} \| \Lambda_k^\frac{1}{4} \tilde{u}_1 \|_{L^2}^\frac{1}{2} \| ( \Lambda_k^\frac{5}{4}, \Lambda_l^\frac{5}{4} ) \tilde{u}_1 \|_{L^2}^\frac{1}{2}. 
	\]
	Use the interpolation inequality
	\[
		\| \Lambda_k^\frac{1}{4} \tilde{u}_1 \|_{L_{x_k}^2} \le C \|\tilde{u}_1 \|_{L_{x_k}^2}^\frac{4}{5} \| \Lambda_k^\frac{5}{4} \tilde{u}_1 \|_{L_{x_k}^2}^\frac{1}{5}, 
	\]
	we obtain that
	\begin{align}
		& | \int \tilde{b}_k\pt_kb_1^{(2)}\tilde{u}_1dx| \le C\|\tilde{b}_k\|_{L^2}\|\tilde{u}_1\|_{L^2}^\frac{2}{5}\|(\Lambda_k^\frac{5}{4},\Lambda_l^\frac{5}{4})\tilde{u}_1\|_{L^2}^\frac{3}{5}\|\Lambda_k^\frac{5}{4}b_1^{(2)}\|_{L^2}^\frac{1}{2}\|\Lambda_k^\frac{5}{4}\pt_1b_1^{(2)}\|_{L^2}^\frac{1}{2} \nonumber \\
		& \le\frac{\nu}{100}\|(\Lambda_k^\frac{5}{4},\Lambda_l^\frac{5}{4})\tilde{u}_1\|_{L^2}^2+C\|\tilde{b}_k\|_{L^2}^\frac{10}{7}\|\tilde{u}_1\|_{L^2}^\frac{4}{7}\|\Lambda_k^\frac{5}{4}b_1^{(2)}\|_{L^2}^\frac{5}{7}\|\Lambda_k^\frac{5}{4}\pt_1b_1^{(2)}\|_{L^2}^\frac{5}{7}. \label{ineq:F2}
	\end{align}
	The first term on the right hand side of \eqref{ineq:F2} is controllable since $i_1\notin\{k,l\}$, and the higher order derivative in the last term is integrable since $j_0 \not= k$. 
	
	{\myemph The case $k\in \mathcal{F}_{3}$. }
	
	Now we have $i_1 = j_0, k=1$. According to $\mathcal{J}$, we must have $i_1 = j_0 = 1$. Due to the solenoidal condition, we have 
	\[
		\int \tilde{b}_1 \pt_1 b_1^{(2)} \tilde{u}_1 dx = - \int \tilde{b}_1 \pt_2 b_2^{(2)} \tilde{u}_1 dx - \int \tilde{b}_1 \pt_3 b_3^{(2)} \tilde{u}_1 dx =: I_a + I_b. 
	\]
	
	Use $f = \tilde{b}_1$, $g = \pt_2 b_2^{(2)}$, $h = \tilde{u}_1$ in Lemma \ref{lem:case:1}, we have 
	\[
		| I_a | \le C \| \tilde{b}_1 \|_{L^2} \| \Lambda_2^\frac{1}{4} \pt_2 b_2^{(2)} \|_{L^2}^\frac{1}{2} \| \Lambda_2^\frac{1}{4} \Lambda_1 \pt_2 b_2^{(2)} \|_{L^2}^\frac{1}{2} \| \Lambda_2^\frac{1}{4} \tilde{u}_1 \|_{L^2}^\frac{1}{2} \| ( \Lambda_2^\frac{5}{4}, \Lambda_3^\frac{5}{4} ) \tilde{u}_1 \|_{L^2}^\frac{1}{2}. 
	\]
	Using the following interpolation inequality
	\[
		\| \Lambda_2^\frac{1}{4} \tilde{u}_1 \|_{L_{x_2}^2} \le C \| \tilde{u}_1 \|_{L_{x_2}^2}^\frac{4}{5} \| \Lambda_2^\frac{5}{4} \tilde{u}_1 \|_{L_{x_2}^2}^\frac{1}{5}, 
	\]
	we obtain that
	\begin{align}
		|I_a| & \le C\|\tilde{b}_1\|_{L^2}\|\tilde{u}_1\|_{L^2}^\frac{2}{5}\|(\Lambda_2^\frac{5}{4},\Lambda_3^\frac{5}{4})\tilde{u}_1\|_{L^2}^\frac{3}{5}\|\Lambda_2^\frac{5}{4}b_2^{(2)}\|_{L^2}^\frac{1}{2}\|\Lambda_2^\frac{5}{4}\pt_1b_2^{(2)}\|_{L^2}^\frac{1}{2} \nonumber \\
		& \le\frac{\nu}{100}\|(\Lambda_2^\frac{5}{4},\Lambda_3^\frac{5}{4})\tilde{u}_1\|_{L^2}^2+C\|\tilde{b}_1\|_{L^2}^\frac{10}{7}\|\tilde{u}_1\|_{L^2}^\frac{4}{7}\|\Lambda_2^\frac{5}{4}b_2^{(2)}\|_{L^2}^\frac{5}{7}\|\Lambda_2^\frac{5}{4}\pt_1b_2^{(2)}\|_{L^2}^\frac{5}{7}. \label{ineq:F3}
	\end{align}
	The first term on the right hand side of \eqref{ineq:F3} is controllable since $i_1\notin\{2,3\}$, and the higher order derivative in the last term is integrable since $j_0\not= 2$. 
	
	The estimates of $I_b$ are similar to those of $I_a$. We omit the details. 
	
	\subsection{The estimates of $L_3$}
	Now we estimate $L_3$. The items in $L_3$ can be expanded as
	\begin{align*}
		L_3 & = \int \tilde{u}_k \partial_k b_m^{(2)} \tilde{b}_m dx 
		= \sum\limits_{k = 1}^{\rm 3} \Big[\int \tilde{u}_k \partial_k b_1^{(2)} \tilde{b}_1 +  \tilde{u}_k \partial_k b_2^{(2)} \tilde{b}_2  + \tilde{u}_k \partial_k b_3^{(2)} \tilde{b}_3 \Big]dx \\
		& =: L_{31} + L_{32} + L_{33}.
	\end{align*}
	We only estimate $L_{31}$ in the following context since $L_{32}$, $L_{33}$ can be argued similarly. Now we separate cases according to the indices $k$. Define
	\[
	\begin{array}{l}
		\mathcal{G}_1 := \{k\in\{1,2,3\} \mid i_k \not= j_0 \}, \\
		\mathcal{G}_2 := \{k\in\{1,2,3\} \mid i_k = j_0 = 1 \}, \\
		\mathcal{G}_3 := \{k\in\{1,2,3\} \mid i_k = j_0 \not= 1\}. 
	\end{array}
	\]
	
	\setcounter{mystep}{0}
	{\myemph The case $k\in \mathcal{G}_1$.}
	
	Since $i_k \not= j_0$, let $l\in\{1,2,3\}\setminus \{ i_k, j_0 \}$, then $\{i_k, j_0, l\}$ is a permutation of $\{1,2,3\}$. Use $f=\tilde{u}_k$, $g=\partial_k b_1^{(2)}$, $h=\tilde{b}_1$ in Lemma \ref{lem:case:4}, we have
	\begin{align}
		& | \int \tilde{u}_k \partial_k b_1^{(2)} \tilde{b}_1 dx| 
		\le C\|\tilde{u}_k\|_{L^2}^\frac{3}{5} \|\Lambda_{j_0}^\frac{5}{4}\tilde{u}_k\|_{L^2}^\frac{2}{5} \|\partial_k b_1^{(2)}\|_{L^2}^\frac{3}{5} \|\Lambda_{i_k}^\frac{5}{4}\partial_k b_1^{(2)}\|_{L^2}^\frac{2}{5} \|\tilde{b}_1\|_{L^2}^\frac{3}{5} \|\Lambda_l^\frac{5}{4}\tilde{b}_1\|_{L^2}^\frac{2}{5} \nonumber \\
		& \le \frac{\nu}{100}\|\Lambda_{j_0}^\frac{5}{4}\tilde{u}_k\|_{L^2}^2 + \frac{\eta}{100}\|\Lambda_l^\frac{5}{4}\tilde{b}_1\|_{L^2}^2 + C\|\tilde{u}_k\|_{L^2} \|\tilde{b}_1\|_{L^2} \|\partial_k b_1^{(2)}\|_{L^2} \|\Lambda_{i_k}^\frac{5}{4}\partial_k b_1^{(2)}\|_{L^2}^\frac{2}{3}. \label{ineq:G1}
	\end{align}
	The first and second terms on the right hand side of \eqref{ineq:G1} are controllable and the higher order derivative in the last term is integrable since $i_k \not= j_0$ and $j_0 \not= l$. 
	
	{\myemph The case $k\in \mathcal{G}_2$.}
	
	From the definition of $\mathcal{J}$, it holds that $k = i_k = j_0 = 1$. Due to the solenoidal condition, we have 
	\[
		- \int \tilde{u}_1\pt_1b_1^{(2)}\tilde{b}_1dx = \int\tilde{u}_1\pt_2b_2^{(2)}\tilde{b}_1dx + \int\tilde{u}_1\pt_3b_3^{(2)}\tilde{b}_1dx =: I_a + I_b.
	\]

	To estimate $I_a$, using $f = \tilde{u}_1$, $g = \pt_1 b_1^{(2)}$, $h = \tilde{b}_1$ in Lemma \ref{lem:case:1}, we have that 
	\begin{align}
		| I_a | & \le C \| \tilde{u}_1 \|_{L^2} \| \Lambda_2^{\frac{1}{4}} \pt_2 b_2^{(2)} \|_{L^2}^{\frac{1}{2}} \| \Lambda_2^{\frac{1}{4}} \Lambda_1 \pt_2 b_2^{(2)} \|_{L^2}^{\frac{1}{2}} \| \Lambda_2^{\frac{1}{4}} \tilde{b}_1 \|_{L^2}^{\frac{1}{2}} \| (\Lambda_2^{\frac{5}{4}}, \Lambda_3^{\frac{5}{4}}) \tilde{b}_1 \|_{L^2}^{\frac{1}{2}} \nonumber \\
		& \le C \| \tilde{u}_1 \|_{L^2} \| \Lambda_2^{\frac{5}{4}} b_2^{(2)} \|_{L^2}^{\frac{1}{2}} \| \Lambda_2^{\frac{5}{4}} \Lambda_1 b_2^{(2)} \|_{L^2}^{\frac{1}{2}} \| \tilde{b}_1 \|_{L^2}^{\frac{2}{5}} \| (\Lambda_2^{\frac{5}{4}}, \Lambda_3^{\frac{5}{4}}) \tilde{b}_1 \|_{L^2}^{\frac{3}{5}} \nonumber \\
		& \le \frac{\eta}{100} \| (\Lambda_2^{\frac{5}{4}}, \Lambda_3^{\frac{5}{4}}) \tilde{b}_1 \|_{L^2}^2 + C \| \tilde{u}_1 \|_{L^2}^{\frac{10}{7}} \| \Lambda_2^{\frac{5}{4}} b_2^{(2)} \|_{L^2}^{\frac{5}{7}} \| \Lambda_2^{\frac{5}{4}} \Lambda_1 b_2^{(2)} \|_{L^2}^{\frac{5}{7}} \| \tilde{b}_1 \|_{L^2}^{\frac{4}{7}}. \label{ineq:G2}
	\end{align}
	The first term on the right hand side of \eqref{ineq:G2} is controllable, and the higher order derivative in the last term is integrable on any finite time interval since $j_0 \not= 2,3$. 
	The estimates of $I_b$ is similar to that of $I_a$. We omit the details for brevity.

	{\myemph The case $k\in \mathcal{G}_3$.}
	
	From the definition of $\mathcal{J}$, it holds now that $k = i_k = j_0 \not= 1$. Let $l\in\{ 1, 2, 3 \} \setminus \{1,k\}$, we have $\{1,k,l\}$ is a permutation of $\{1,2,3\}$. Use $f=\tilde{u}_k$, $g=\partial_k b_1^{(2)}$, $h=\tilde{b}_1$ in Lemma \ref{lem:case:4}, we have
	\begin{align}
		& | \int \tilde{u}_k \partial_k b_1^{(2)} \tilde{b}_1 dx|
		\le C\|\tilde{u}_k\|_{L^2}^\frac{3}{5} \|\Lambda_1^\frac{5}{4}\tilde{u}_k\|_{L^2}^\frac{2}{5} \|\partial_k b_1^{(2)}\|_{L^2}^\frac{3}{5} \|\Lambda_k^\frac{5}{4}\partial_k b_1^{(2)}\|_{L^2}^\frac{2}{5} \|\tilde{b}_1\|_{L^2}^\frac{3}{5} \|\Lambda_l^\frac{5}{4}\tilde{b}_1\|_{L^2}^\frac{2}{5} \nonumber \\
		& \le \frac{\nu}{100}\|\Lambda_1^\frac{5}{4}\tilde{u}_k\|_{L^2}^2 + \frac{\eta}{100}\|\Lambda_l^\frac{5}{4}\tilde{b}_1\|_{L^2}^2 + C\|\tilde{u}_k\|_{L^2} \|\tilde{b}_1\|_{L^2} \|\partial_k b_1^{(2)}\|_{L^2} \|\Lambda_k^\frac{5}{4}\partial_k b_1^{(2)}\|_{L^2}^\frac{2}{3}. \label{ineq:G3}
	\end{align}
	The first and second terms on the right hand side of \eqref{ineq:G3} are controllable since $i_k \not= 1$ and $j_0 \not= l$. The last term $\|\Lambda_k^\frac{5}{4}\partial_k b_1^{(2)}\|$ is bounded according to the assumption in Theorem \ref{thm:uniqueness:MHD}.
	
	\subsection{The estimates of $L_4$}
	Now we estimate $L_4$. The items in $L_4$ can be expanded as
	\begin{align*}
		L_4 & = \int \tilde{b}_k \partial_k u_m^{(2)} \tilde{b}_m dx
		=\sum\limits_{k = 1}^{\rm 3} \Big[\int \tilde{b}_k \partial_k u_1^{(2)} \tilde{b}_1 + \tilde{b}_k \partial_k u_2^{(2)} \tilde{b}_2 + \tilde{b}_k \partial_k u_3^{(2)} \tilde{b}_3\Big] dx \\
		& =: L_{41} + L_{42} + L_{43}.
	\end{align*}
	We only estimate $L_{41}$ in the following context since $L_{42}$, $L_{43}$ can be argued similarly. Now we separate cases according to the indices $k$. Define
	\[
	\begin{array}{l}
		\mathcal{H}_1 := \{k\in\{1,2,3\} \mid i_1 \not= j_0 \}, \\
		\mathcal{H}_2 := \{k\in\{1,2,3\} \mid i_1 = j_0, k \not= 1 \}, \\
		\mathcal{H}_3 := \{k\in\{1,2,3\} \mid i_1 = j_0, k = 1 \}. 
	\end{array}
	\]
	
	\setcounter{mystep}{0}
	{\myemph The case $k\in \mathcal{H}_{1}$. }
	
	Now we have $j_0 \not= i_1$. Let $l\in\{ 1, 2, 3 \} \setminus \{ i_1, j_0 \}$, then $\{ i_1, j_0, l \}$ is a permutation of $\{ 1, 2, 3 \}$. Use $f=\tilde{b}_k$, $g = \partial_k u_1^{(2)}$, $h = \tilde{b}_1$ in Lemma \ref{lem:case:4}, we have
	\begin{align}
		& | \int \tilde{b}_k \partial_k u_1^{(2)} \tilde{b}_1 dx| 
		\le C\|\tilde{b}_k\|_{L^2}^\frac{3}{5} \|\Lambda_l^\frac{5}{4}\tilde{b}_k\|_{L^2}^\frac{2}{5} \|\partial_k u_1^{(2)}\|_{L^2}^\frac{3}{5} \|\Lambda_{j_0}^\frac{5}{4}\partial_k u_1^{(2)}\|_{L^2}^\frac{2}{5} \|\tilde{b}_1\|_{L^2}^\frac{3}{5} \|\Lambda_{i_1}^\frac{5}{4}\tilde{b}_1\|_{L^2}^\frac{2}{5} \nonumber \\
		& \le \frac{\eta}{100}\|\Lambda_l^\frac{5}{4}\tilde{b}_k\|_{L^2}^2 + \frac{\eta}{100}\|\Lambda_{i_1}^\frac{5}{4}\tilde{b}_1\|_{L^2}^2 + C\|\tilde{b}_k\|_{L^2} \|\tilde{b}_1\|_{L^2} \|\partial_k u_1^{(2)}\|_{L^2} \|\Lambda_{j_0}^\frac{5}{4}\partial_k u_1^{(2)}\|_{L^2}^\frac{2}{3}. \label{ineq:H1}
	\end{align}
	The first and second terms on the right hand side of \eqref{ineq:H1} are controllable, and the higher order derivative in the last term is integrable on any finite time interval since $j_0\not\in\{l, i_1\}$. 
	
	{\myemph The case $k\in \mathcal{H}_{2}$. }
	
	Now we have $i_1 = j_0, k \not= 1$. According to the set $\mathcal{J}$, we must have $i_1 = j_0 = 1$. Let $l\in\{1,2,3\}\setminus \{1,k\}$, then $\{1,k,l\}$ is a permutation of $\{1,2,3\}$. Use $f=\tilde{b}_k$, $g=\pt_ku_1^{(2)}$, $h=\tilde{b}_1$ in Lemma \ref{lem:case:1}, we have 
	\[
		| \int \tilde{b}_k \pt_k u_1^{(2)} \tilde{b}_1 dx | \le C \| \tilde{b}_k \|_{L^2} \| \Lambda_k^\frac{1}{4} \pt_k u_1^{(2)} \|_{L^2}^\frac{1}{2} \| \Lambda_k^\frac{1}{4} \Lambda_1 \pt_k u_1^{(2)} \|_{L^2}^\frac{1}{2} \| \Lambda_k^\frac{1}{4} \tilde{b}_1 \|_{L^2}^\frac{1}{2} \| ( \Lambda_k^\frac{5}{4}, \Lambda_l^\frac{5}{4} ) \tilde{b}_1 \|_{L^2}^\frac{1}{2}. 
	\]
	Using the interpolation inequality
	\[
		\| \Lambda_k^\frac{1}{4} \tilde{b}_1 \|_{L_{x_k}^2} \le C \| \tilde{b}_1 \|_{L_{x_k}^2}^\frac{4}{5} \| \Lambda_k^\frac{5}{4} \tilde{b}_1 \|_{L_{x_k}^2}^\frac{1}{5}, 
	\]
	we obtain that 
	\begin{align}
		& |\int \tilde{b}_k \pt_k u_1^{(2)} \tilde{b}_1 dx | \le C \| \tilde{b}_k \|_{L^2} \| \tilde{b}_1 \|_{L^2}^\frac{2}{5} \| ( \Lambda_k^\frac{5}{4}, \Lambda_l^\frac{5}{4}) \tilde{b}_1 \|_{L^2}^\frac{3}{5} \| \Lambda_k^\frac{5}{4} u_1^{(2)} \|_{L^2}^\frac{1}{2} \| \Lambda_k^\frac{5}{4} \pt_1 u_1^{(2)} \|_{L^2}^\frac{1}{2} \nonumber \\
		& \le \frac{\eta}{100} \| ( \Lambda_k^\frac{5}{4}, \Lambda_l^\frac{5}{4} ) \tilde{b}_1 \|_{L^2}^2 + C \| \tilde{b}_k \|_{L^2}^\frac{10}{7} \| \tilde{b}_1 \|_{L^2}^\frac{4}{7} \| \Lambda_k^\frac{5}{4} u_1^{(2)} \|_{L^2}^\frac{5}{7} \| \Lambda_k^\frac{5}{4} \pt_1 u_1^{(2)} \|_{L^2}^\frac{5}{7}. \label{ineq:H2}
	\end{align}
	The first term on the right hand side of \eqref{ineq:H2} is controllable, and the higher order derivative in the last term is integrable on any finite time interval since $j_0 = i_1 \notin \{ k, l \}$. 
	
	{\myemph The case $k\in \mathcal{H}_{3}$. }
	
	Now we have $i_1 = j_0, k=1$. According to the set $\mathcal{J}$, we must have $i_1 = j_0 = 1$. Due to the solenoidal condition, we have 
	\[
		\int \tilde{b}_1 \pt_1 u_1^{(2)} \tilde{b}_1 dx = - \int \tilde{b}_1 \pt_2 u_2^{(2)} \tilde{b}_1 dx - \int \tilde{b}_1 \pt_3 u_3^{(2)} \tilde{b}_1 dx =: I_a + I_b. 
	\]
	
	From the definition of $\mathcal{J}$, we must have $i_2 \not= 1$. Let $l\in\{1,2,3\}\setminus \{i_1,i_2\}$, we have $\{i_1,i_2,l\}$ is a permutation of $\{1,2,3\}$. Use $f=\tilde{b}_1$, $g=\partial_2 u_2^{(2)}$, $h=\tilde{b}_1$ in Lemma \ref{lem:case:4}, we have
	\begin{align}
		|I_a| & \le C\|\tilde{b}_1\|_{L^2}^\frac{3}{5} \|\Lambda_{i_2}^\frac{5}{4}\tilde{b}_1\|_{L^2}^\frac{2}{5} \|\partial_2 u_2^{(2)}\|_{L^2}^\frac{3}{5} \|\Lambda_{i_1}^\frac{5}{4}\partial_2 u_2^{(2)}\|_{L^2}^\frac{2}{5} \|\tilde{b}_1\|_{L^2}^\frac{3}{5} \|\Lambda_l^\frac{5}{4}\tilde{b}_1\|_{L^2}^\frac{2}{5} \nonumber \\
		& \le \frac{\eta}{100}\|\Lambda_{i_2}^\frac{5}{4}\tilde{b}_1\|_{L^2}^2 + \frac{\eta}{100}\|\Lambda_l^\frac{5}{4}\tilde{b}_1\|_{L^2}^2 + C\|\tilde{b}_1\|_{L^2}^2 \|\partial_2 u_2^{(2)}\|_{L^2} \|\Lambda_{i_1}^\frac{5}{4}\partial_2 u_2^{(2)}\|_{L^2}^\frac{2}{3}. \label{ineq:H3}
	\end{align}
	The first and second terms on the right hand side of \eqref{ineq:H3} are controllable since $j_0 \notin \{i_2,l\}$. The higher order derivative in the last term is integrable on any finite time interval since $i_2 \not= i_1$. 
	
	Note that $i_3 \not= 1$, then the estimate of $I_b$ are similar to those of $I_a$. The details are omitted for brevity.
	
	\vspace{10pt}
	
	Combining all the above estimates for $L_1$ through $L_4$ together, we obtain that
	\begin{align*}
		& \frac{d}{dt}\| (\widetilde{\bm{u}}, \widetilde{\bm{b}}) \|_{L^2}^2 + \nu \sum_{j\not=i_1}\| \Lambda_j ^\frac{5}{4} \tilde{u}_1 \|_{L^2}^2 + \nu \sum_{j\not=i_2}\| \Lambda_j ^\frac{5}{4} \tilde{u}_2 \|_{L^2}^2 + \nu \sum_{j\not=i_3}\| \Lambda_j ^\frac{5}{4} \tilde{u}_3 \|_{L^2}^2 \\
		& + \eta \sum_{k=1}^{3} \sum_{j\not=j_0} \| \Lambda_j ^\frac{5}{4} \tilde{b}_k \|_{L^2}^2 \le B(t) \| (\widetilde{\bm{u}}, \widetilde{\bm{b}}) \|_{L^2}^2,
	\end{align*}
	where
	\begin{align*}
		B(t) & = C \Big(\|\nabla {\bm u}^{(2)}\|_{L^2}(\sum_{j\not=i_1}\| \Lambda_j ^\frac{5}{4} \nabla u_1^{(2)} \|_{L^2}^\frac{2}{3}+\sum_{j\not=i_2}\| \Lambda_j ^\frac{5}{4} \nabla u_2^{(2)} \|_{L^2}^\frac{2}{3}+\sum_{j\not=i_3}\| \Lambda_j ^\frac{5}{4} \nabla u_3^{(2)} \|_{L^2}^\frac{2}{3}) \\
		& + \|\nabla {\bm b}^{(2)}\|_{L^2} \sum_{k=1}^{3} \sum_{j \not= j_0 } \| \Lambda_j ^\frac{5}{4} \nabla b_1^{(2)} \|_{L^2}^\frac{2}{3} + \sum_{j\not=i_1}\| \Lambda_j ^\frac{5}{4} u_1^{(2)} \|_{L^2}^\frac{5}{7}\| \Lambda_j ^\frac{5}{4} \nabla u_1^{(2)} \|_{L^2}^\frac{5}{7} \\
		& + \sum_{j\not=i_2}\| \Lambda_j ^\frac{5}{4} u_2^{(2)} \|_{L^2}^\frac{5}{7}\| \Lambda_j ^\frac{5}{4} \nabla u_2^{(2)} \|_{L^2}^\frac{5}{7} + \sum_{j\not=i_3}\| \Lambda_j ^\frac{5}{4} u_3^{(2)} \|_{L^2}^\frac{5}{7}\| \Lambda_j ^\frac{5}{4} \nabla u_3^{(2)}\|_{L^2}^\frac{5}{7} \\
		& + \sum_{k=1}^{3} \sum_{j\not=j_0} \| \Lambda_j ^\frac{5}{4} b_1^{(2)} \|_{L^2}^\frac{5}{7}\| \Lambda_j ^\frac{5}{4} \nabla b_1^{(2)} \|_{L^2}^\frac{5}{7} + \sum_{k\not=1}\delta_{ki_1}\delta_{ki_k}\|\nabla u_1^{(2)} \|_{L^2} \|\Lambda_k^\frac{9}{4} u_1^{(2)} \|_{L^2}^\frac{2}{3} \\
		& + \sum_{k\not=2}\delta_{ki_2}\delta_{ki_k}\|\nabla u_2^{(2)} \|_{L^2} \|\Lambda_k^\frac{9}{4} u_2^{(2)} \|_{L^2}^\frac{2}{3} + \sum_{k\not=3}\delta_{ki_3}\delta_{ki_k}\|\nabla u_3^{(2)} \|_{L^2} \|\Lambda_k^\frac{9}{4} u_3^{(2)} \|_{L^2}^\frac{2}{3} \\
		& + \sum_{j=1}^3 \sum_{k\not=j}\delta_{kj_0}\delta_{ki_k}\|\nabla b_j^{(2)}\|_{L^2} \|\Lambda_k^\frac{9}{4} b_j^{(2)}\|_{L^2}^\frac{2}{3}
		\Big).
	\end{align*}
	The regularity assumption on ${\bm u}^{(2)}$ and ${\bm b}^{(2)}$ in \eqref{eq:reg}-\eqref{eq:uniqueness:cond:2} ensures the integrability of $B(t)$ for any fixed $T>0$, the Gr\"{o}nwall's inequality then implies that $ ({\bm u}^{(1)}, {\bm b}^{(1)}) \equiv ({\bm u}^{(2)}, {\bm b}^{(2)})$ if $({\bm u}_0^{(1)}, {\bm b}_0^{(1)}) \equiv ({\bm u}_0^{(2)}, {\bm b}_0^{(2)})$. This completes the proof of uniqueness.
\hfill$\square$

Corollary \ref{cor:MHD} can be readily derived through Theorem \ref{thm:existence:MHD} and Theorem \ref{thm:uniqueness:MHD}. 


\section*{Appendix}
\renewcommand{\thesection}{A}

This section first introduces several fundamental inequalities required, including: Sobolev embedding inequalities, Gagliardo-Nirenberg inequalities and Minkowski inequalities. The proof of the main theorems in this paper relies on estimates of a family of trilinear terms. We summarize the estimates into four types, which are presented in Lemma \ref{lem:case:1} through Lemma \ref{lem:case:4}.
\begin{lemma}[Sobolev embedding]\label{lem:Sobolev:ineq}
	Suppose $2 \le p \le \infty $ and $s > n( \frac{1}{2} - \frac{1}{p} )$. Then there exists a constant $C$, which depends only on $n$, $p$ and $s$, such that
	\[
	\|f\|_{L^p(\mathbb{R}^n)} \le C 
	\|f\|_{L^2(\mathbb{R}^n)}^{ 1 - \frac{n}{s} (\frac{1}{2} -\frac{1}{p})} \cdot \|\Lambda^s f \|_{L^2 (\mathbb{R}^n)}^{\frac{n}{s}( \frac{1}{2} - \frac{1}{p})}, \qquad f \in {H^s}( \mathbb{R}^n ).
	\]
	Moreover, if $2\le p < +\infty$, the above inequality also holds for $s = n (\frac{1}{2} - \frac{1}{p})$. 
\end{lemma}

The following four lemmas will be utilized repeatedly in the proof of our main results.

\begin{lemma}\label{lem:case:1}
	Suppose that functions $f$, $\Lambda_i^{\frac{1}{4}}g$, $\Lambda_i^{\frac{1}{4}} \Lambda_j g$, $\Lambda_i^{\frac{1}{4}} h$ and $(\Lambda_i^{\frac{5}{4}},\Lambda_k^{\frac{5}{4}}) h$ belong to $L^2(\mathbb{R}^3)$ for some $(i,j,k)$ being a permutation of $\{1,2,3\}$. Then 
	\[
	\Big| \int_{\mathbb{R}^3} f g h dx \Big| \le \| f \|_{L^2} \cdot \| \Lambda_i^{\frac{1}{4}} g \|_{L^2}^{\frac{1}{2}} \cdot \| \Lambda_i^{\frac{1}{4}} \Lambda_j g \|_{L^2}^{\frac{1}{2}} \cdot \| \Lambda_i^{\frac{1}{4}} h \|_{L^2}^{\frac{1}{2}} \cdot \| (\Lambda_i^{\frac{5}{4}},\Lambda_k^{\frac{5}{4}}) h \|_{L^2}^{\frac{1}{2}}.
	\]
\end{lemma}
\begin{proof}
	By H{\"o}lder's inequality, we have 
	\[
	\Big|\int_{\mathbb{R}^3} f g h dx \Big| \le \| f \|_{L^2} \| g \|_{L_{x_k}^2 L_{x_j}^\infty L_{x_i}^4} \| h \|_{L_{x_k}^\infty L_{x_j}^2 L_{x_i}^4}. 
	\]
	By Lemma \ref{lem:Sobolev:ineq} and Minkowski's inequality, we have
	\[
	\begin{array}{l}
		\| g \|_{L_{x_k}^2 L_{x_j}^\infty L_{x_i}^4} \le C \| \Lambda_i^\frac{1}{4} g \|_{L_{x_k x_i}^2 L_{x_j}^\infty}
		\le C \| \Lambda_i^\frac{1}{4} g \|_{L^2}^\frac{1}{2} \| \Lambda_i^\frac{1}{4} \Lambda_j g\|_{L^2}^\frac{1}{2}, \\
		\| h \|_{L_{x_k}^\infty L_{x_j}^2 L_{x_i}^4} \le C \| \Lambda_i^\frac{1}{4} h \|_{L_{x_j x_i}^2 L_{x_k}^\infty} 
		\le C\| \Lambda_i^\frac{1}{4} h \|_{L^2}^\frac{1}{2}\| ( \Lambda_i^\frac{5}{4}, \Lambda_k^\frac{5}{4}) h \|_{L^2}^\frac{1}{2}.
	\end{array}
	\]
	This completes the proof of Lemma \ref{lem:case:1}.
\end{proof}
\noindent
\emph{Remark}: It is sometimes important that we keep $\|\Lambda_i^{\frac{1}{4}} \Lambda_j g\|_{L^2}$ in the above lemma, although it is obvious that $\|\Lambda_i^{\frac{1}{4}} \Lambda_j g\|_{L^2} \le C \|(\Lambda_i^{\frac{5}{4}},\Lambda_j^{\frac{5}{4}})g\|_{L^2}$. 
\begin{lemma}\label{lem:case:2}
	Suppose that functions $\Lambda_i^{\frac{1}{4}}f$, $\Lambda_j^{\frac{1}{4}}g$, $\Lambda_j^{\frac{1}{4}}\Lambda_kg$, $h $ and $(\pt_i,\pt_j)h$ belong to $L^2(\mathbb{R}^3)$ for some $(i,j,k)$ being a permutation of $\{1,2,3\}$. Then
	\[
	\Big|\int_{\mathbb{R}^3} f g h dx\Big|\le\|\Lambda_i^{\frac{1}{4}}f\|_{L^2}\cdot\|\Lambda_j^{\frac{1}{4}}g\|_{L^2}^{\frac{1}{2}}\cdot\|\Lambda_j^{\frac{1}{4}}\Lambda_kg\|_{L^2}^{\frac{1}{2}}\cdot\|h\|_{L^2}^{\frac{1}{2}}\cdot\|(\pt_i,\pt_j)h\|_{L^2}^{\frac{1}{2}}.
	\]
\end{lemma}
\begin{proof}
	By H{\"o}lder's inequality, we have 
	\[
		\Big|\int_{\mathbb{R}^3} f g h dx\Big|\le \|f\|_{{L_{x_k x_j}^2}{L_{x_i}^4}} \|g\|_{L_{x_k}^\infty L_{x_j}^4 L_{x_i}^2 } \|h\|_{L_{x_k}^2 L_{x_j x_i}^4}. 
	\]
	By Lemma \ref{lem:Sobolev:ineq}, Gagliardo-Nirenberg's inequality and Minkowski's inequality, we have
	\[
	\begin{array}{l}
		\| f \|_{{L_{x_j x_k}^2}{L_{x_i}^4}} \le C \big{\|} \| \Lambda_i^\frac{1}{4}f\|_{L_{x_i}^2} \big{\|}_{L_{x_j x_k}^2}
		= C\| \Lambda_i^\frac{1}{4}f\|_{L^2}, \\
		\| g \|_{L_{x_i}^2 L_{x_k}^\infty L_{x_j}^4} \le C\| \Lambda_j^\frac{1}{4}g\|_{L_{x_i x_j}^2 L_{x_k}^\infty}
		\le C\| \Lambda_j^\frac{1}{4}g\|_{L^2}^\frac{1}{2}\| \Lambda_j^\frac{1}{4} \Lambda_k g\|_{L^2}^\frac{1}{2}, \\
		\| h \|_{L_{x_k}^2 L_{x_i x_j}^4} \le C {\big{\|} \| h\|_{L_{x_ix_j}^2}^\frac{1}{2}
			\| (\pt_i,\pt_j)h\|_{L_{x_i x_j}^2}^\frac{1}{2} \big{\|}}_{L_{x_k}^2} \le C \| h\|_{L^2}^\frac{1}{2}
		\| (\pt_i,\pt_j)h\|_{L^2}^\frac{1}{2}.
	\end{array}
	\]
	This completes the proof of Lemma \ref{lem:case:2}.
\end{proof}
\begin{lemma}\label{lem:case:3}
	Suppose that functions $\Lambda_i^{\frac{1}{4}}f$, $g$, $\Lambda_k^{\frac{5}{4}}g$, $ h$ and $(\Lambda_i^{\frac{5}{4}},\Lambda_j^{\frac{5}{4}}) h$ belong to $L^2(\mathbb{R}^3)$ for some $(i,j,k)$ being a permutation of $\{1,2,3\}$. Then
	\[
		\Big| \int_{\mathbb{R}^3} f g h dx \Big| \le \| \Lambda_i^{\frac{1}{4}} f \|_{L^2} \cdot \| g \|_{L^2}^{\frac{3}{5}} \cdot \| \Lambda_k^{\frac{5}{4}} g \|_{L^2}^{\frac{2}{5}} \cdot \| h \|_{L^2}^{\frac{2}{5}} \cdot \| ( \Lambda_i^{\frac{5}{4}}, \Lambda_j^{\frac{5}{4}}) h \|_{L^2}^{\frac{3}{5}}.
	\]
\end{lemma}
\begin{proof}
	By H{\"o}lder's inequality, we have 
	\[
		\Big|\int_{\mathbb{R}^3} f g h dx\Big|\le \|f\|_{{L_{x_k x_j}^2}{L_{x_i}^4}} \|g\|_{L_{x_k}^\infty L_{x_j x_i}^2} \|h\|_{L_{x_k}^2 L_{x_j}^\infty L_{x_i}^4}. 
	\]
	By Lemma \ref{lem:Sobolev:ineq} and Minkowski's inequality, we have
	\[
	\begin{array}{l}
		\| f \|_{{L_{x_j x_k}^2}{L_{x_i}^4}} \le C\| \Lambda_i^\frac{1}{4}f\|_{L^2}, \\
		\| g \|_{L_{x_k}^\infty L_{x_j x_i}^2 } 
		\le C {\big{\|} \| g\|_{L_{x_k}^2}^\frac{3}{5}
			\| \Lambda_k^{\frac{5}{4}} g \|_{L_{x_k}^2}^{\frac{2}{5}} \big{\|}}_{L_{x_j x_i}^2} \le C \| g \|_{L^2}^{\frac{3}{5}}
		\| \Lambda_k^{\frac{5}{4}} g \|_{L^2}^{\frac{2}{5}}.
	\end{array}
	\]
	Similarly, by Lemma \ref{lem:Sobolev:ineq}, Gagliardo-Nirenberg's inequality, Minkowski's inequality and the interpolation inequality
	\[
		\|\Lambda_i^\frac{1}{4} h\|_{L^2} \le C \|h\|_{L^2}^\frac{4}{5} \|\Lambda_i^\frac{5}{4} h\|_{L^2}^\frac{1}{5}, 
	\]
	we have that 
	\begin{align*}
		\| h \|_{L_{x_k}^2 L_{x_j}^\infty L_{x_i}^4} & \le C \| \Lambda_i^\frac{1}{4} h \|_{L_{x_i x_k}^2 L_{x_j}^\infty}
		\le C \| \Lambda_i^\frac{1}{4} h \|_{L^2}^\frac{1}{2} \| \Lambda_i^{\frac{1}{4}} \Lambda_j h \|_{L^2}^\frac{1}{2} \\
		& \le C \| h \|_{L^2}^\frac{2}{5}\| \Lambda_i^\frac{5}{4} h\|_{L^2}^\frac{1}{10}\|(\Lambda_i^{\frac{5}{4}},\Lambda_j^{\frac{5}{4}}) h\|_{L^2}^\frac{1}{2}
		\le C \| h \|_{L^2}^\frac{2}{5} \| ( \Lambda_i^{\frac{5}{4}}, \Lambda_j^{\frac{5}{4}}) h \|_{L^2}^\frac{3}{5}, 
	\end{align*}
	This completes the proof of Lemma \ref{lem:case:3}.
\end{proof}	

\begin{lemma}\label{lem:case:4}
	Suppose that functions $f$, $\Lambda_i^\frac{5}{4} f$, $g$, $\Lambda_j^{\frac{5}{4}}g$, $h$ and $\Lambda_k^{\frac{5}{4}}h$ belong to $L^2(\mathbb{R}^3)$ for some $(i,j,k)$ being a permutation of $\{1,2,3\}$. Then 
	\[
		\Big| \int_{\mathbb{R}^3} f g h dx \Big| \le \| f \|_{L^2}^\frac{3}{5} \cdot \| \Lambda_i^\frac{5}{4} f \|_{L^2}^\frac{2}{5} \cdot \| g \|_{L^2}^\frac{3}{5} \cdot \| \Lambda_j^\frac{5}{4} g \|_{L^2}^\frac{2}{5} \cdot \| h \|_{L^2}^\frac{3}{5} \cdot \| \Lambda_k^\frac{5}{4} h \|_{L^2}^\frac{2}{5}.
	\]
\end{lemma}
\begin{proof}
	By H{\"o}lder's inequality, we have 
	\[
		\Big|\int_{\mathbb{R}^3} f g h dx\Big|\le\|f\|_{L_{x_i}^\infty L_{x_jx_k}^2} \|g\|_{L_{x_i}^2 L_{x_j}^\infty L_{x_k}^2} \|h\|_{L_{x_ix_j}^2 L_{x_k}^\infty }. 
	\]
	By Lemma \ref{lem:Sobolev:ineq} and Minkowski's inequality, we have
	\[
		\|f\|_{L_{x_i}^\infty L_{x_j x_k}^2 } 
		\le C \big{\|} \| f\|_{L_{x_i}^2}^\frac{3}{5}
		\| \Lambda_i^{\frac{5}{4}}f\|_{L_{x_i}^2}^\frac{2}{5} \big{\|}_{L_{x_j x_k}^2}\le C \| f\|_{L^2}^\frac{3}{5}
		\| \Lambda_i^{\frac{5}{4}}f\|_{L^2}^\frac{2}{5}.
	\]
	Similarly, we have that
	\[
		\| g \|_{ L_{x_i x_k}^2 L_{x_j}^\infty} 
		\le C \| g \|_{L^2}^\frac{3}{5} \| \Lambda_j^{\frac{5}{4}} g \|_{L^2}^\frac{2}{5}, \qquad
		\| h \|_{ L_{x_i x_j}^2 L_{x_k}^\infty} 
		\le C \| h \|_{L^2}^\frac{3}{5} \| \Lambda_k^{\frac{5}{4}} h \|_{L^2}^\frac{2}{5}.
	\]
	This completes the proof of Lemma \ref{lem:case:4}.
\end{proof}

\bibliographystyle{plain}
\bibliography{partial_diff}

\end{document}